\documentclass{article}
\usepackage{graphicx}
\usepackage{placeins}
\usepackage[margin=1.3in]{geometry}
\usepackage{amsmath}
\usepackage{amsthm}
\usepackage{amssymb}
\usepackage{mathtools}
\usepackage{enumitem}
\usepackage{booktabs}
\usepackage{multirow}
\usepackage{empheq}
\usepackage{parskip}
\usepackage{algorithm}
\usepackage{algpseudocodex}
\usepackage{hyperref}
\usepackage{caption}
\usepackage{subcaption}
\usepackage{siunitx}
\usepackage{pdflscape} 
\usepackage{longtable}
\usepackage{xcolor}
\usepackage{colortbl}
\usepackage{adjustbox}
\usepackage{xspace}

\newcommand{\R}{\mathbb{R}}

\theoremstyle{definition}
\newtheorem{theorem}{Theorem}[section]
\newtheorem{lemma}[theorem]{Lemma}
\newtheorem{proposition}[theorem]{Proposition}
\newtheorem{corollary}[theorem]{Corollary}
\newtheorem{definition}[theorem]{Definition}

\newcommand{\cpot}{\textsc{cpo3}\xspace}

\title{Strongly Monotone Central Path Equations for Quadratic Cone Programs}
\author{%
  Michael Cummins\textsuperscript{\textasteriskcentered,\textdagger}
  \and
  Matteo Santamaria\textsuperscript{\textasteriskcentered,\textdaggerdbl}
  \and
  Brendan O'Donoghue\textsuperscript{\textsection}
}
\date{} 

\begin{document}

\maketitle
{\renewcommand{\thefootnote}{\fnsymbol{footnote}}%
\footnotetext[1]{These authors contributed equally to this work.}%
\footnotetext[2]{Imperial College London, UK. \texttt{m.cummins24@imperial.ac.uk}}%
\footnotetext[3]{University of California, Berkeley, USA. \texttt{matteosantama@gmail.com}}%
\footnotetext[4]{Google DeepMind, London, UK. \texttt{bodonoghue85@gmail.com}}}

\begin{abstract}
We introduce a canonical central path for the homogeneous embedding of quadratic cone programs, defined by strongly monotone equations whose modulus is the barrier parameter. We prove that the path exists uniquely, is smooth, and lies on a sphere of fixed radius. Under mild regularity assumptions, its limit points recover optimal solutions or certificates of infeasibility. Exploiting the fact that the embedding operator has a closed-form resolvent, we derive Josephy--Newton steps that retain the nonlinear embedding equations and linearize only the complementarity conditions. Each step reduces to quasi-definite linear system solves and a scalar quadratic. Using this, we derive a new predictor--corrector interior-point method with a surrogate formulation that allows the predictor and corrector to share a single sparse factorization. We provide an open-source implementation, \cpot{}, for linear and quadratic programs. On 675 instances from eight benchmark collections, with presolve and crossover disabled, \cpot{} achieves the highest verified solve count (631) and the lowest penalized shifted geometric mean runtime among the tested open-source and commercial solvers. It is the only tested solver to solve all 138 Maros--M\'esz\'aros instances to high accuracy. Moreover, we establish Euclidean and local-norm error bounds that motivate a warm-start heuristic for the proposed algorithm, yielding up to 2.3$\times$ median runtime speedup over cold starts on perturbed Netlib and Maros--M\'esz\'aros problems. 
\end{abstract}

\section{Introduction}

In this paper we consider the standard-form Quadratic Cone Program (QCP) and its dual,

\begin{equation*}
\begin{minipage}{0.48\linewidth}
\begin{equation}
\begin{aligned}
    \min_{x, s} &\quad \frac{1}{2} x^\top P x + c^\top x, \\
    \text{s.t.} &\quad A x + s = b, \\ &\quad s \in \mathcal{K},
\end{aligned}
\label{eq:primal_qp}
\end{equation}
\end{minipage}\hfill
\begin{minipage}{0.48\linewidth}
\begin{equation}
\begin{aligned}
    \max_{\xi, y} &\quad -\frac{1}{2} \xi^\top P \xi - b^\top y, \\
    \text{s.t.} &\quad P \xi + A^\top y + c = 0, \\ &\quad y \in \mathcal{K}^*.
\end{aligned}
\label{eq:dual_qp}
\end{equation}
\end{minipage}
\end{equation*}
where $x \in \R^n$ and $s \in \R^m$ are the primal variables, and
$y \in \R^m$ and $\xi \in \R^n$ are the dual variables.
We assume strong duality  holds between \eqref{eq:primal_qp} and \eqref{eq:dual_qp}, and their minimizers are attained if both are feasible.
The matrix $P \in \R^{n \times n}$
is assumed to be symmetric positive semidefinite, and we write
$\mathcal K^*$ for the dual of the closed convex cone
$\mathcal K \subseteq \R^m$. We represent $z \in \{0,\ldots,m\}$
equality constraints by a zero-cone factor $\{0\}^z$ in $\mathcal K$,
whose dual is $\R^z$. The remaining factor, when present, is assumed
to be closed, convex and pointed, with nonempty interior; we refer
to it as the proper cone factor. This formulation includes linear
and quadratic programs as special cases and, through conic
reformulation, a broad class of convex problems arising in control \cite{RawlingsMayneDiehl2017,BorrelliBemporadMorari2017},
signal processing \cite{combettes1996convex}, finance \cite{markowitz1952portfolio,boyd2017multiperiod} and machine learning \cite{cortes1995support}.

Interior-point methods (IPMs) have long been a leading approach
for solving QCPs, owing to their strong theoretical guarantees
\cite{nesterov1994,nesterov1998primal,skajaa2015homogeneous} and
numerical robustness \cite{wright1995stability,wright1998ill}.
Path-following IPMs are \emph{homotopy methods}: they seek solutions
of the Karush--Kuhn--Tucker (KKT) optimality conditions by approximately
following a \emph{central path} of perturbed KKT systems as the barrier
parameter decreases to zero \cite{wright1997primal,ye1997interior}. Standard methods
compute search directions by linearizing these equations. The
resulting linear system is a primary source of numerical difficulty:
its coefficient matrix is typically symmetric indefinite and can
become severely ill-conditioned as the iterates approach a solution.
One remedy is symmetric indefinite factorization with numerical
pivoting \cite{golub2013matrix}, which improves stability but can
change the sparsity pattern of the factors, increasing 
solve time. To avoid numerical pivoting, some implementations
explicitly regularize the KKT system
\cite{domahidi2013ecos,goulart2024,chari2025,zanetti2025factorisation},
allowing a symbolic factorization to be reused across iterations.
When regularization is introduced solely as a numerical perturbation,
the computed direction need not solve the original KKT system, and
iterative refinement is often used to reduce the residual in that
system. In practice, the speed and robustness of a solver depend
strongly on how these systems are handled. In this work, we derive 
strongly monotone central-path equations that yield `regularized' quasi-definite matrices by construction. This admits pivot-free factorization in exact
arithmetic with a symmetric permutation chosen once at the start
of the solve \cite{vanderbei1995}.

A second challenge is detecting infeasibility of
\eqref{eq:primal_qp} or \eqref{eq:dual_qp}. A standard mechanism is
the homogeneous embedding \cite{andersen1999homogeneous}, which
augments the problem so that optimal solutions and infeasibility
certificates can be recovered from a single system of equations.
For QCPs, the embedding operator is monotone but generally nonlinear \cite{goulart2024,andersen1999homogeneous}.
Existing embedding-based IPMs typically either linearize the
embedding \cite{goulart2024} or reformulate the quadratic objective
using an additional second-order cone constraint
\cite{domahidi2013ecos}. Although both approaches are effective,
we exploit additional structure in the embedding to obtain a more
practical and efficient interior-point method.

Our contributions are as follows.
\begin{itemize}
    \item We define a new, \emph{canonical} central path for the
    homogeneous embedding of \eqref{eq:primal_qp} defined by strongly monotone equations. We establish
    existence and uniqueness for every positive barrier parameter
    and show that the path is smooth and lies on a sphere of fixed
    radius. When a finite primal--dual solution or a strict
    infeasibility certificate exists, its limit points recover
    optimal solutions or certificates of primal or dual
    infeasibility.

    \item We derive a modified Mehrotra predictor--corrector method
    that computes step directions by linearizing only the
    complementarity conditions and solving the remaining nonlinear
    equations in closed form. Each step reduces to quasi-definite
    linear solves and a scalar quadratic equation, while a surrogate
    formulation allows the predictor and corrector to share a
    single sparse factorization.

    \item We provide an open-source Rust implementation of our IPM
    for linear and quadratic programs. On 675 instances drawn from
    eight benchmark collections, with presolve and crossover
    disabled, \cpot{} achieves the highest verified solve count
    and the lowest penalized shifted geometric mean runtime among
    the tested open-source and commercial solvers.

    \item We use strong monotonicity and barrier geometry to derive
    Euclidean and local-norm error bounds for the distance to the
    central path. These bounds motivate a warm-start heuristic
    that yields up to 2.3$\times$ median runtime speedup on perturbed Maros--M\'esz\'aros and Netlib problems.

\end{itemize}
We call this path canonical because, for fixed problem data and barrier,
it is unique and independent of the initial point.
The central-path analysis applies to any proper cone equipped with
a logarithmically homogeneous self-concordant barrier, with equality
constraints handled separately. The algorithm and implementation
are presented for linear and quadratic programs; their extension
to general cones is left to future work. 

\paragraph{Structure:}Section~\ref{sec:preliminaries} introduces mathematical preliminaries for our work. The
homogeneous embedding of \eqref{eq:primal_qp} is introduced in Section~\ref{sec:homogenous}, and we define its
central path in Section~\ref{sec:central_path_section}. The new path complicates the direct application of the standard Mehrotra-style predictor--corrector algorithm; Section~\ref{sec:ipm_section} shows how to adapt the method to overcome this difficulty. We describe the implementation of our solver in Section~\ref{sec:implementation} and benchmark against commercial and open source solvers in Section~\ref{sec:numerical}.

\paragraph{Notation:} We denote the set of nonnegative and strictly positive real numbers by $\R_{\geq 0}$ and $\R_{>0}$, respectively.
A tuple $(a, b, c)$ denotes the stacked column vector $[a^\top \; b^\top \; c^\top]^\top$. We write $x \perp y$ if two vectors $x, y \in \mathbb{R}^d$ satisfy $x^\top y = 0$. For
some set $S$, we define its interior and relative interior by $\operatorname{int}(S)$ and $\operatorname{relint}(S)$,
respectively. For a set-valued operator $T$ defined on $\R^d$, the domain is denoted by
$\operatorname{dom} T = \{x \,:\, T(x) \neq \emptyset\}$. The zero set of an operator $T$ is denoted by $\operatorname{zer}(T) \coloneq \{x \,:\, 0 \in T(x)\}$. We denote the Euclidean norm of a vector $x \in \R^d$ by
$\| x \|$. For some matrix $M \succ 0$, we define the norm $\|x\|_M = \sqrt{x^\top M x}$. The space of $n$-times
continuously differentiable functions with domain $S$ is denoted by $C^n(S)$.

\subsection{Related work}
\paragraph{Regularized Interior-Point Methods:}
Vanderbei's quasi-definiteness result \cite{vanderbei1995} motivated
many IPM variants that regularize the KKT system statically or dynamically
\cite{saunders1996solving,altman1999regularized,friedlander2012primal,pougkakiotis2019dynamic}.
For methods admitting a proximal-point interpretation
\cite{rockafellar1976monotone,rockafellar1976augmented}, the
regularized direction is itself an algorithmic step and need not
be refined to solve the unregularized KKT system. The regularization
parameters nevertheless require suitable selection and update rules.
Regularization proportional to the complementarity parameter also appears
in central paths for nonlinear complementarity problems
\cite{zhao2006path} and linear programming \cite{castro2010existence}.
Closely related is the work of Pougkakiotis and Gondzio
\cite{pougkakiotis2021interior,pougkakiotis2022interior}, who developed
an Interior-Point Proximal Method of Multipliers (IP-PMM). Inspired by Armand and Omheni \cite{armand2017mixed}, they solve \eqref{eq:primal_qp} via the proximal method of multipliers, and compute approximate solutions to the bound constrained subproblems using an IPM.

\paragraph{Homogeneous Algorithms:}
The homogeneous self-dual embedding was introduced for LPs
in \cite{ye1994homogeneous}. Homogeneous algorithms were subsequently
developed for monotone complementarity problems
\cite{andersen1999homogeneous} and complementarity problems over
symmetric cones \cite{yoshise2007interior}.
For linear cone programs ($P=0$), Nesterov, Todd and Ye
\cite{nesterov1999infeasible} defined a central path for the
homogeneous self-dual embedding and established complexity
guarantees for a path-following IPM. Skajaa and Ye
\cite{skajaa2015homogeneous} extended this approach to nonsymmetric
cones, with more practical variants later developed for various
cones \cite{dahl2022primal,chen2025efficient,he2024qics}.
The homogeneous QCP embedding considered here retains the quadratic
objective but is generally nonlinear and no longer self-dual. O'Donoghue
\cite{odonoghue2021} showed that the embedding operator admits a
maximal monotone extension and used operator splitting to solve
the resulting system, extending the original SCS algorithm from
linear to quadratic cone programs \cite{ocpb:16}.
An alternative homogeneous QP formulation with infeasibility detection
\cite{raghunathan2021homogeneous} assumes a positive-definite reduced
Hessian and supports an IPM with warm starts
\cite{raghunathan2022homogeneous}.

\paragraph{Solvers:}
IPM-based solvers for \eqref{eq:primal_qp} differ in their use of
homogeneous embeddings and their treatment of the step equations.
PIQP \cite{schwan2023piqp}, HiGHS \cite{huangfu2018highs} and
QOCO \cite{chari2025} handle quadratic objectives without conic
reformulation, but their IPMs do not use homogeneous embeddings.
PIQP adapts the proximal method of
\cite{pougkakiotis2021interior}, while HiGHS includes the HiPO
interior-point solver \cite{zanetti2025factorisation}, which
dynamically regularizes the KKT system in the spirit of
\cite{altman1999regularized}.
QOCO follows an approach similar to
\cite[Section 5]{vandenberghe2010cvxopt}, supports second-order
cone constraints with a quadratic objective, and uses
regularization and iterative refinement.
Tulip \cite{tanneau2021design} regularizes a homogeneous LP embedding
around the current iterate.
MOSEK \cite{andersen2000mosek,andersen2003implementing} and
ECOS \cite{domahidi2013ecos} use homogeneous embeddings for conic
formulations with linear objectives, handling quadratic objectives
through conic reformulation.
Clarabel \cite{goulart2024} supports a range of cones while retaining
the quadratic objective through the embedding described in
\cite{odonoghue2021}. Its step computation linearizes both the
embedding and the complementarity conditions, and solves the
resulting KKT system using regularization and iterative refinement.
Gurobi \cite{gurobi2025} offers a homogeneous variant of its barrier
method for diagnosing infeasibility and unboundedness, while
CPLEX \cite{cplex2025} provides a primal--dual barrier solver.

\section{Preliminaries}\label{sec:preliminaries}
This section presents mathematical preliminaries fundamental to both interior point methods and the understanding of our central path.

\subsection{Infeasibility certificates}\label{sec:certificates}

The primal or dual problem is \emph{infeasible} if its constraints
admit no feasible point. In the strongly infeasible case, this can
be certified by a strongly separating hyperplane. Define the sets
of primal and dual infeasibility certificates,
\begin{gather*}
    \mathcal{P} \coloneq \Big\{ y \in \R^m :
    A^\top y = 0, \:\: y \in \mathcal{K}^*, \:\: b^\top y < 0 \Big\}, \\
    \mathcal{D} \coloneq \Big\{ x \in \R^n :
    Px = 0, \:\: -Ax \in \mathcal{K}, \:\: c^\top x < 0 \Big\}.
\end{gather*}
Any $y \in \mathcal P$ certifies strong infeasibility of
\eqref{eq:primal_qp}, while any $x \in \mathcal D$ certifies
strong infeasibility of \eqref{eq:dual_qp} \cite{odonoghue2021}. In other words,
the existence of one of these certificates precludes the existence of a primal-dual
solution to \eqref{eq:primal_qp} and \eqref{eq:dual_qp}, and vice versa.

\subsection{Barrier functions}
Barrier functions are used in interior point methods to convert conic constraints into equations, enabling the definition of a central path as a nonlinear system of equations.

\begin{definition}\label{def:v-lhscb}
A three-times continuously differentiable, self-concordant \cite[\S 2.1]{nesterov1994} function
$f : \operatorname{int}(\mathcal{K}) \to \R$ is a
\emph{logarithmically homogeneous self-concordant barrier of degree $\nu \ge 1$}
($\nu$-LHSCB) for a proper cone $\mathcal{K}$ if $f(x_k) \to \infty$ for every sequence
$\{x_k\} \subset \operatorname{int}(\mathcal{K})$ converging to the boundary of $\mathcal{K}$, and
\begin{equation}\label{eq:lhscb_lh}
    f(\lambda x) = f(x) - \nu \log \lambda,
    \quad \forall x \in \operatorname{int}(\mathcal{K}),\ \lambda > 0.
\end{equation}
\end{definition}
Self-concordance implies $\nabla^2 f(x)\succ0$ and logarithmic homogeneity implies
$x^\top\nabla f(x)=-\nu$. The conjugate barrier
\begin{equation}\label{eq:conj_def}
    f^*(y)\coloneq\sup_{x\in\operatorname{int}(\mathcal K)}
    \{-y^\top x-f(x)\},
\end{equation}
is a $\nu$-LHSCB for $\mathcal K^*$ \cite{nesterov1997self}.
For every $y\in\operatorname{int}(\mathcal K^*)$, the supremum is
attained uniquely at $-\nabla f^*(y)\in\operatorname{int}(\mathcal K)$.
In our setting, we apply the preceding definition to the proper cone factor and extend the conjugate
barrier constantly over the unrestricted equality multipliers if present.
Its gradient and Hessian therefore have zero equality components.
If there is no proper cone factor, we use the convention
$f=f^*=0$ and $\nu=0$.

\subsection{Monotone Operators}
The central path analysis in Section~\ref{sec:central_path_section} relies on the machinery of monotone operators, which we make explicit in this section.
A set-valued operator $T$ on $\R^d$ assigns to each $u\in\R^d$ a
(possibly empty) subset $T(u)\subseteq\R^d$. Its graph is
$\operatorname{gra}T\coloneq\{(u,v):v\in T(u)\}$,
a subset of $\R^d\times\R^d$. We say $T$ is single-valued if
$T(u)$ is a singleton for every $u\in\operatorname{dom}T$.
An operator $T$ is monotone on $\R^d$ if
\begin{equation}\label{eq:monotone_ineq}
    (u-\bar u)^\top(v-\bar v)\geq0,
    \qquad
    v\in T(u),\quad \bar v\in T(\bar u).
\end{equation}
For a single-valued $T$, condition \eqref{eq:monotone_ineq} reads
$(u-\bar u)^\top\bigl(T(u)-T(\bar u)\bigr)\geq0$
for all $u,\bar u\in\operatorname{dom}T$. Furthermore, $T$ is
$\mu$-strongly monotone, for some $\mu>0$, if
\begin{equation}\label{eq:strongly_monotone_ineq}
    (u-\bar u)^\top(v-\bar v)\geq\mu\|u-\bar u\|^2,
    \qquad
    v\in T(u),\quad \bar v\in T(\bar u).
\end{equation}
An operator $T$ is maximal monotone on $\R^d$ if it is monotone
and its graph is not strictly contained in the graph of another
monotone operator. That is, if $T$ is not maximal monotone,
there must exist $(u,v)\notin\operatorname{gra}T$ such that
$\operatorname{gra}T\cup\{(u,v)\}$ is the graph of a monotone
operator. If an operator $\mathcal T$ is maximal monotone and
$\operatorname{gra}T\subseteq\operatorname{gra}\mathcal T$,
then $\mathcal T$ is a maximal monotone extension of $T$.
The interested reader is directed to
\cite{Bauschke_Combettes_2017} and \cite{ryu2016primer} for treatments of this topic.

\section{Homogeneous Embedding}\label{sec:homogenous}

After eliminating $\xi$, the KKT optimality conditions for \eqref{eq:primal_qp} are \cite{boyd2004convex}
\begin{subequations}
    \label{eq:kkt_conditions}
    \begin{equation}
        \begin{bmatrix}
            P & A^\top \\
            -A & 0
        \end{bmatrix}
        \begin{bmatrix}
            x \\
            y
        \end{bmatrix}
        +
        \begin{bmatrix}
            c \\
            b
        \end{bmatrix}
        =
        \begin{bmatrix}
            0 \\
            s
        \end{bmatrix},
    \end{equation}
    \begin{equation}
        s \perp y, \qquad
        s \in \mathcal{K}, \qquad
        y \in \mathcal{K}^*.
    \end{equation}
\end{subequations}
Any triple $(x,y,s)$ satisfying \eqref{eq:kkt_conditions}
yields a primal optimal solution $(x,s)$ and a dual optimal
solution $(\xi,y)=(x,y)$. If either problem is infeasible,
however, the KKT system has no solution.

The homogeneous embedding augments this system with two
nonnegative scalars, $\tau$ and $\kappa$, so that optimal
solutions and certificates of strong infeasibility can be
represented within a single system. Following
\cite{andersen1999homogeneous,odonoghue2021}, we first define
the single-valued operator
$Q : \R^n \times \R^m \times \R_{>0} \to \R^{n+m+1}$ by
\begin{equation}\label{eq:Q_operator}
    Q(x, y, \tau) \coloneq \begin{bmatrix}
        \phantom{+}P\phantom{^\top} & \phantom{+}A^\top & \phantom{+}c \\
        -A\phantom{^\top} & \phantom{+}0\phantom{^\top} & \phantom{+}b \\
        -c^\top & -b^\top & \phantom{+}0
    \end{bmatrix} \begin{bmatrix}
        x \\
        y \\
        \tau
    \end{bmatrix} - \frac{1}{\tau} \begin{bmatrix}
        0 \\
        0 \\
        x^\top P x
    \end{bmatrix}.
\end{equation}
The operator is positively homogeneous of degree one and is \emph{pointwise orthogonal}, \emph{i.e.}, satisfies
\begin{equation}\label{eq:conservative}
    (x,y,\tau)^\top Q(x,y,\tau)=0.
\end{equation}
It is also monotone on its domain
\cite[Lemma 4.1]{odonoghue2021}, but not maximal monotone.
In particular, its domain excludes $\tau=0$, where the embedding
must represent infeasibility certificates. To include these
boundary points, we use the following extension.

\begin{corollary}[{\cite[Lemma 4.4]{odonoghue2021}}]
    \label{cor:unique_Q}
    The operator $Q$ admits a unique maximal monotone extension
    $\mathcal Q$ whose domain is contained in
    $\R^n \times \R^m \times \R_{\geq0}$.
\end{corollary}

This extension agrees with $Q$ for $\tau>0$. Its graph is
$\operatorname{gra}\mathcal Q
=\operatorname{gra}Q\cup\operatorname{gra}\mathcal I$,
where the boundary component $\mathcal I$ is given by
\begin{equation}\label{eq:I_operator}
    \mathcal I(x,y,0)
    =
    \left\{
    (A^\top y, -Ax, \eta)
        : \eta \leq -c^\top x-b^\top y
    \right\},
\end{equation}
with domain $\bigl\{(x,y,0)\in\R^{n+m+1}:Px=0\bigr\}$.
Outside its domain, $\mathcal I$ is empty. The complete
homogeneous embedding is then
\begin{equation}
    \label{eq:homo_embedding}
    \begin{bmatrix}
        0 \\
        s \\
        \kappa
    \end{bmatrix}
    \in \mathcal Q(x,y,\tau), \quad
    (s,\tau)\perp(y,\kappa), \quad
    s\in\mathcal K, \quad y\in\mathcal K^*, \quad
    \tau\geq0, \quad \kappa\geq0.
\end{equation}
Note that $\mathcal{Q}$ is pointwise orthogonal, which implies the second requirement, but we include it here for completeness.
The embedding always has the trivial zero solution so is always feasible.
Solutions with $\tau+\kappa>0$ yield either optimal solutions
or infeasibility certificates, as the following result shows.

\begin{lemma}[Solution or certificate recovery]
    \label{lem:asymp_feasible}
    Let $(x,y,s,\tau,\kappa)$ solve
    \eqref{eq:homo_embedding}.
    \begin{enumerate}[label=\textit{(\roman*)}]
        \item If $\tau>0$, then $\kappa=0$, and $(x/\tau,s/\tau)$ and
        $(x/\tau,y/\tau)$ solve
        \eqref{eq:primal_qp} and \eqref{eq:dual_qp},
        respectively.

        \item If $\kappa>0$, then $\tau=0$, and
        at least one of $y\in\mathcal P$ or
        $x\in\mathcal D$ holds.
    \end{enumerate}
\end{lemma}

\begin{proof}
    Cone membership and complementarity imply
    $s^\top y=\tau\kappa=0$.
\begin{enumerate}[label=\textit{(\roman*)}]
    \item If $\tau>0$, then $\kappa=0$ and
    $\mathcal Q(x,y,\tau)
    =\{Q(x,y,\tau)\}$.
    Dividing the first two block equations by $\tau$
    and rescaling the complementarity condition shows that
    $(x/\tau,y/\tau,s/\tau)$ satisfies the KKT conditions
    \eqref{eq:kkt_conditions}.
    \item If $\kappa>0$, then $\tau=0$, $Px=0$, and $\mathcal Q(x,y,\tau)
    =\mathcal I(x,y,\tau)$. This implies that $A^\top y=0$, $-Ax=s\in\mathcal K$, $y\in\mathcal K^*$, and
    \[
        c^\top x+b^\top y
        \leq-\kappa<0.
    \]
    Hence at least one of $b^\top y$ or $c^\top x$
    is negative, yielding a certificate in $\mathcal P$
    or $\mathcal D$, respectively.
\end{enumerate}

\end{proof}

The recovery result does not classify solutions with
$\tau=\kappa=0$, including the trivial zero solution.
Section~\ref{sec:limiting} establishes conditions that exclude
this case for limit points of the central path.

The central path and the accepted iterates of our proposed IPM maintain
$\tau>0$, so we formulate the central-path equations using
the single-valued map $Q$. The extension $\mathcal Q$ gives
meaning to boundary solutions at $\tau=0$ and is used in the
theoretical analysis of Section~\ref{sec:central_path_section}.
A detailed derivation of the extension is given in
\cite[Sections 4.3--4.5]{odonoghue2021}.

\section{Central path}\label{sec:central_path_section}

For brevity, let
\[
    u\coloneq(x,y,\tau), \qquad v\coloneq(0,s,\kappa),
\]
and define the composite cone
$\mathcal C\coloneq\R^n\times\mathcal K^*\times\R_{\geq0}$
and its dual
$\mathcal C^*\coloneq\{0\}^n\times\mathcal K\times\R_{\geq0}$.
In this notation, the homogeneous embedding
\eqref{eq:homo_embedding} is
\begin{equation}\label{eq:embedding_compact}
    v\in\mathcal Q(u),\quad
    u\perp v,\quad
    u\in\mathcal C,\quad
    v\in\mathcal C^*.
\end{equation}
Let $f$ be a $\nu$-LHSCB for the proper cone factor of
$\mathcal K$. As in Section~\ref{sec:preliminaries}, we extend
its conjugate barrier $f^*$ constantly over the unrestricted
equality multipliers, so the corresponding gradient and Hessian
components are zero. Equality rows do not contribute to $\nu$;
if there is no proper cone factor, we use $f^*=0$ and $\nu=0$.
Define $\phi:\operatorname{int}(\mathcal C)\to\R$ by
\begin{equation}\label{eq:phi-def}
    \phi(u)\coloneq f^*(y)-\log\tau.
\end{equation}
For each barrier parameter $\mu>0$, we define the central-path
equations on $u\in\operatorname{int}(\mathcal C)$ by
\begin{subequations}\label{eq:central_path}
    \begin{equation}\label{eq:Q_resolvent}
        Q(u)+\mu u=v,
    \end{equation}
    \begin{equation}\label{eq:ys_self-concordance}
        v=-\mu\nabla\phi(u).
    \end{equation}
\end{subequations}
Since $\tau>0$ on this domain, the single-valued operator $Q$
suffices to define the path. The same parameter $\mu$
regularizes the embedding equations and perturbs the
complementarity conditions.
Equation~\eqref{eq:ys_self-concordance} ensures
$v\in\operatorname{relint}(\mathcal C^*)$ and gives $s=-\mu\nabla f^*(y)$, $\tau\kappa=\mu$, and  $s^\top y=\mu\nu$.
On a symmetric proper cone factor with its canonical logarithmic
barrier we have,
\begin{equation}\label{eq:product_complementarity}
    s=-\mu\nabla f^*(y)\quad \Leftrightarrow\quad s\circ y=\mu e,
\end{equation}
where $\circ$ and $e$ denote the Jordan product and identity
on the proper cone factor \cite{faraut1994}, extended by zero
on equality rows.
For the nonnegative orthant, the Jordan product is componentwise,
giving $s_i y_i=\mu$ on each inequality row. This is the
specialization used in Section~\ref{sec:ipm_section}.
Equality slacks remain zero, while equality multipliers receive
regularization through $\mu I$ but no barrier term.

Substituting \eqref{eq:ys_self-concordance} into
\eqref{eq:Q_resolvent} gives a description in terms of $u$ alone:
\begin{equation}\label{eq:central_path_reduced}
    T_\mu(u)\coloneq Q(u)+\mu u+\mu\nabla\phi(u)=0.
\end{equation}
This is the main equation that defines our central path.
Since $Q$ and $\nabla\phi$ are monotone, $T_\mu$ is
$\mu$-strongly monotone:
\[
    (u-\bar u)^\top\bigl(T_\mu(u)-T_\mu(\bar u)\bigr)
    \geq\mu\|u-\bar u\|^2,
    \qquad
    u,\bar u\in\operatorname{int}(\mathcal C).
\]
This property underpins the uniqueness and error-bound results
below. The embedding also retains its computational structure:
the scaled resolvent $(\mathcal Q+\mu I)^{-1}$ admits a
closed-form evaluation \cite{odonoghue2021}.
Section~\ref{sec:ipm_section} exploits this structure to
linearize only the $s \circ y$ complementarity conditions while
retaining the nonlinear embedding equations.
The regularization also fixes the scale of the central-path
points, as the following identity shows.
\begin{proposition}[Sphere identity]
    \label{prop:sphere}
    Let $\mu>0$. If $u\in\operatorname{int}(\mathcal C)$
    satisfies $T_\mu(u)=0$, then
    \[
        \|u\|^2=\nu+1.
    \]
\end{proposition}
\begin{proof}
    Logarithmic homogeneity of $\phi$ gives
    $u^\top\nabla\phi(u)=-(\nu+1)$, while
    $u^\top Q(u)=0$.
    Taking the inner product of $T_\mu(u)=0$ with $u$ yields
    \[
        0
        =u^\top Q(u)+\mu\|u\|^2
         +\mu u^\top\nabla\phi(u)
        =\mu\bigl(\|u\|^2-(\nu+1)\bigr).
        \qedhere
    \]
\end{proof}
Strong monotonicity ensures that $T_\mu$ has at most one zero.
The following result establishes existence for every $\mu>0$,
without requiring a KKT point or an infeasibility certificate.

\begin{theorem}[Existence and uniqueness of the central path]
    \label{thm:central_path_existence}
    For every $\mu>0$, the equation $T_\mu(u)=0$ admits a
    unique solution
    $u_\mu \coloneq (x_\mu,y_\mu,\tau_\mu)\in\operatorname{int}(\mathcal C)$.
\end{theorem}
The central path is also smooth, with its differentiability
inherited from the barrier.

\begin{theorem}[Smoothness of the central path]
    \label{thm:central-path-smoothness}
    Suppose
    $f^*\in C^{d+1}(\operatorname{int}(\mathcal K^*))$
    for some integer $d\geq2$. Then the central-path solution
    map $\mu\mapsto u_\mu$ defined by
    \eqref{eq:central_path_reduced} is of class $C^d$ on
    $\R_{>0}$.
\end{theorem}

The proofs of these two theorems are given in
Appendix~\ref{app:path}. For the standard symmetric-cone
barriers, $f^*\in C^\infty(\operatorname{int}(\mathcal K^*))$,
so $(\mu\mapsto u_\mu)\in C^\infty(\R_{>0})$.
We next derive error bounds that motivate warm starting,
before examining the limits as $\mu\downarrow0$ and
$\mu\to\infty$.

\subsection{Error Bounds}\label{subsec:error_bounds}

For fixed problem data and a chosen barrier, our central
path is independent of the initial point. This differs from
homogeneous central-path constructions anchored at an initial
point \cite{nesterov1998primal,skajaa2015homogeneous}.
The following residual-based bounds quantify proximity to
the path point at a prescribed parameter $\mu$ and motivate
the warm-start heuristic described in
Section~\ref{subsec:warm-start}.

\begin{proposition}[Euclidean error bound]
    \label{prop:error_bound}
    Let $\mu>0$ and let $u_\mu$ be the unique solution
    of \eqref{eq:central_path_reduced}. Then, for every
    $u\in\operatorname{int}(\mathcal C)$,
    \begin{equation}\label{eq:error-bound}
        \|u-u_\mu\|
        \leq\frac{1}{\mu}\|T_\mu(u)\|.
    \end{equation}
\end{proposition}
\begin{proof}
    If $u=u_\mu$, the bound is immediate.
    Otherwise, strong monotonicity,
    $T_\mu(u_\mu)=0$, and Cauchy--Schwarz give
    \[
        \mu\|u-u_\mu\|^2
        \leq
        (u-u_\mu)^\top T_\mu(u)
        \leq
        \|u-u_\mu\|\,\|T_\mu(u)\|.
    \]
\end{proof}

The factor $1/\mu$ can make the Euclidean bound conservative
for small $\mu$, and Euclidean distance does not reflect
the local barrier geometry. To incorporate this geometry,
define the self-concordant, strongly convex function
\[
    \Psi(u)\coloneq\tfrac12\|u\|^2+\phi(u),
    \qquad
    \nabla^2\Psi(u)=I+\nabla^2\phi(u)\succeq I.
\]
Although $\phi$ is flat in the $x$ and equality-multiplier
directions, the quadratic term makes $\nabla^2\Psi(u)$
positive definite on the full space. Its Hessian therefore
defines a local norm, with the inverse Hessian defining
the corresponding dual norm.

\begin{theorem}[Local-norm error bound]
    \label{thm:barrier_error_bound}
    Let $\mu>0$ and $u\in\operatorname{int}(\mathcal C)$,
    and define
    \begin{equation}\label{eq:eta_def}
        \eta_\mu(u)
        \coloneq
        \frac{1}{\mu}
        \|T_\mu(u)\|_{\left( \nabla^2\Psi(u) \right)^{-1}}.
    \end{equation}
    If $\eta_\mu(u)<1$, then
    \[
        \|u-u_\mu\|
        \leq
        \|u-u_\mu\|_{\nabla^2\Psi(u)}
        \leq
        \frac{\eta_\mu(u)}{1-\eta_\mu(u)}.
    \]
\end{theorem}

The proof is given in Appendix~\ref{app:path}.
The residual measure $\eta_\mu(u)$ thus provides a sufficient
condition for proximity to the central path in the local norm.
We use this measure to motivate the warm-start heuristic in
Section~\ref{subsec:warm-start}. 

\subsection{Limiting Behavior}\label{sec:limiting}

We now relate the central path to solutions of the homogeneous
embedding as $\mu\downarrow0$.
Write $v_\mu=(0,s_\mu,\kappa_\mu)
\coloneq-\mu\nabla\phi(u_\mu)$ for the corresponding slack vector using \eqref{eq:ys_self-concordance}.
Proposition~\ref{prop:sphere} bounds $u_\mu$ uniformly;
the following result also bounds $v_\mu$ as $\mu\downarrow0$
and identifies the limiting system.

\begin{lemma}[Limit points of the central path]
    \label{lem:limit_points}
    For every $\mu_0>0$, the family
    $\{(u_\mu,v_\mu):0<\mu\leq\mu_0\}$ is bounded.
    Consequently, the central path has at least one limit
    point as $\mu\downarrow0$, and every such limit point
    solves the homogeneous embedding
    \eqref{eq:embedding_compact}.
\end{lemma}

The proof is given in Appendix~\ref{app:path}.
This result shows that the central path has at least one limit point as
$\mu\downarrow0$, and that every such limit point solves the homogeneous
embedding~\eqref{eq:embedding_compact}. To recover a KKT solution or an
infeasibility certificate, we must also exclude limit points with
$\tau=\kappa=0$. Theorem~\ref{thm:nondegenerate_limits} does so under strong duality and the existence of a strong infeasibility certificate.

\begin{theorem}[Limiting solutions]
    \label{thm:nondegenerate_limits}
    Let $\mu>0$ and let $\tau_\mu$ and $\kappa_\mu$ be the corresponding values on the central path.
    \begin{enumerate}[label=\textit{(\roman*)}]
        \item If $(x,y,s)$ is a KKT
        point, then
        \[
            \tau_\mu
            \geq
            \frac{1}{
                \sqrt{\nu+1}
                \sqrt{1+\|x\|^2+\|y\|^2}
            }.
        \]

        \item If $y\in\mathcal P$, then
        \[
            \tau_\mu
            \leq
            \frac{\mu\sqrt{\nu+1}\|y\|}{-b^\top y},
            \qquad
            \kappa_\mu
            \geq
            \frac{-b^\top y}{\sqrt{\nu+1}\|y\|}.
        \]

        \item If $x\in\mathcal D$, then
        \[
            \tau_\mu
            \leq
            \frac{\mu\sqrt{\nu+1}\|x\|}{-c^\top x},
            \qquad
            \kappa_\mu
            \geq
            \frac{-c^\top x}{\sqrt{\nu+1}\|x\|}.
        \]
    \end{enumerate}
\end{theorem}

The proof is given in Appendix~\ref{app:path}.
In case \textit{(i)}, every limit point has $\tau>0$
and therefore recovers a primal--dual optimal solution by
Lemma~\ref{lem:asymp_feasible}. In cases \textit{(ii)} and
\textit{(iii)}, every limit point has $\kappa>0$ and
therefore recovers an infeasibility certificate.

Finally, we consider the opposite limit, $\mu\to\infty$.
The Euclidean error bound identifies a
\emph{large-parameter endpoint} that depends only on the
chosen cone and barrier, not on the problem data
$(P,A,b,c)$. We use this endpoint as the default initialization
in our implementation.

\begin{proposition}[The large-parameter endpoint]
    \label{prop:large_mu_endpoint}
    Let $u_\infty$ be the unique minimizer of $\Psi$ over
    $\operatorname{int}(\mathcal C)$. Then
    \[
        u_\infty=(0,y_\infty,1),\qquad
        y_\infty+\nabla f^*(y_\infty)=0,\qquad
        \|u_\infty\|^2=\nu+1.
    \]
    In particular, the equality-multiplier component of
    $y_\infty$ is zero. For every $\mu>0$,
    \begin{equation}\label{eq:large_mu_bound}
        \|u_\mu-u_\infty\|
        \leq\frac{1}{\mu}\|Q(u_\infty)\|.
    \end{equation}
    Consequently, $u_\mu\to u_\infty$ as
    $\mu\to\infty$.
    For $\mathcal K=\{0\}^z\times\R_{\geq0}^{m-z}$
    with the standard logarithmic barrier,
    $y_\infty=(0_z,\mathbf1_{m-z})$.
\end{proposition}
\begin{proof}
    Extending $\Psi$ by $+\infty$ outside
    $\operatorname{int}(\mathcal C)$ gives a proper,
    closed, strongly convex function. It therefore has
    a unique minimizer in its domain.
    The stationarity condition
    \[
        \nabla\Psi(u_\infty)
        =u_\infty+\nabla\phi(u_\infty)=0
    \]
    gives the stated form of $u_\infty$, including zero
    equality multipliers because the corresponding
    gradient components of $f^*$ vanish.
    Taking its inner product with $u_\infty$ yields
    the sphere identity.
    Since $T_\mu(u_\infty)=Q(u_\infty)$,
    \eqref{eq:large_mu_bound} follows from
    Proposition~\ref{prop:error_bound}.
\end{proof}

\begin{figure}[!htbp]
\centering
\includegraphics[width=\linewidth]{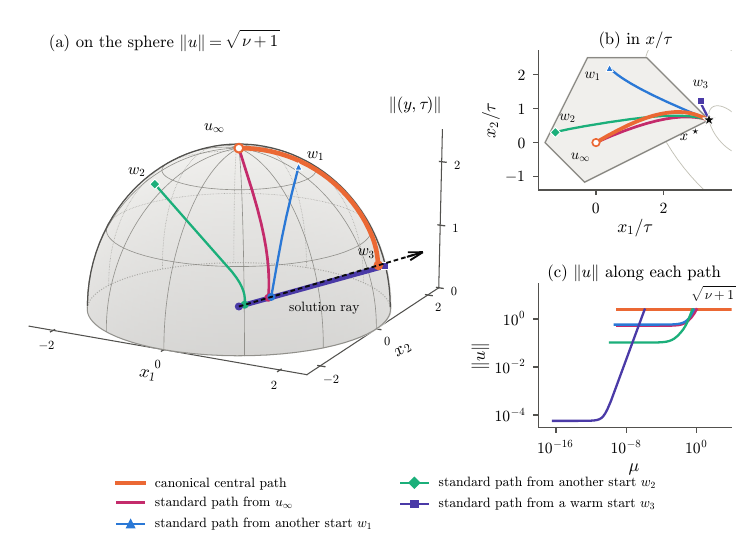}
\caption{
Canonical (orange) and standard homogeneous central paths for a QP with
two variables and five inequality constraints.
(a)~Paths in coordinates $(x_1,x_2,\|(y,\tau)\|)$, where the normalization
sphere appears as a hemisphere.
(b)~The corresponding primal coordinates $x/\tau$.
(c)~The norm $\|u\|$ as a function of $\mu$.}
\label{fig:central_paths}
\end{figure}

Together with Theorem~\ref{thm:nondegenerate_limits} and Proposition~\ref{prop:sphere}, this shows that the
central path joins a point fixed by the cone and barrier to its limit points
as $\mu\downarrow0$, and that it lies on the sphere $\|u\|^2=\nu+1$
throughout. Figure~\ref{fig:central_paths} illustrates this for a QP with two
variables, where we plot our canonical central path and compare it to three standard central paths from different starting points \cite{andersen1999homogeneous,goulart2024}. The canonical path starts at $u_\infty$ and meets the solution ray where the ray leaves the sphere. The comparison paths are anchored at their starting points, leave the
sphere, and meet the ray at points that depend on the starting point. In this example, the warm-started standard path approaches
the solution ray near the origin, so recovering $x/\tau$ to a given accuracy
requires a much smaller $\mu$. Our path and its scale are the same for
every starting point.
All starting points are normalized to $\|u_0\|^2=\nu+1$, and each curve
is traced until $x/\tau$ is within $10^{-6}$ of $x^\star$.
The warm-started standard path reaches this accuracy only at
$\mu\approx10^{-16}$.

\section{Interior-point method}\label{sec:ipm_section}

In this section, we develop a primal--dual predictor--corrector
interior-point method for following the central path described in
Section~\ref{sec:central_path_section}. To simplify the notation
and implementation discussion, we specialize to
$\mathcal K=\{0\}^{z}\times\R_{\geq0}^{m-z}$ with the standard
logarithmic barrier on the inequality rows. Extensions to other
symmetric cones can be developed using standard Nesterov-Todd scaling matrices
\cite{nesterov1998primal}. 
For (much) more detail on interior point methods see, \emph{e.g.}, \cite{wright1997primal,ye1997interior,nesterov1994,renegar2001mathematical,roos2005interior}.

Standard predictor--corrector methods compute step directions
by linearizing the central-path equations. Here, although
\eqref{eq:Q_operator} is nonlinear, we can retain the embedding
equations and linearize only the complementarity conditions
$s\circ y=\mu e$. The resulting step can be computed in closed
form and interpreted as a \emph{Josephy--Newton} iteration
\cite{josephy1979newton}.

\subsection{Computing Step Directions}\label{subsec:josephy_newton}

For an equation $F(u,s)+G(u,s)=0$, the Josephy--Newton method
computes a trial point $(u_{k+1},s_{k+1})$ from starting point $(u_{k},s_{k})$ by solving the partially linearized equation
\begin{equation}\label{eq:josephy_newton}
    F(u_k,s_k)+\nabla F(u_k,s_k)
    \begin{bmatrix}
        u_{k+1}-u_k\\
        s_{k+1}-s_k
    \end{bmatrix}
    +G(u_{k+1},s_{k+1})=0,
\end{equation}
giving the step directions $\Delta u_k=u_{k+1}-u_k$ and
$\Delta s_k=s_{k+1}-s_k$.
Related Newton methods for self-concordant inclusions retain the monotone
operator while linearizing the barrier gradient \cite{trandinh2019self}.
We apply this framework to \eqref{eq:central_path}, treating
the primal slack $s$ as an independent variable alongside
$u=(x,y,\tau)$. Equality slacks remain fixed at zero and
equality multipliers are unrestricted. We use componentwise
products and inverses on inequality rows, with $s\circ y$, $e$
and $y^{-1}$ set to zero on equality rows.
Throughout this subsection, the subscript $k+1$ denotes an
undamped trial point.

Eq.~\eqref{eq:ys_self-concordance} is equivalent to $s\circ y=\mu e$ and $\tau \kappa = \mu$. For a target $\mu>0$, substituting $\kappa=\mu/\tau$ into
\eqref{eq:Q_resolvent} eliminates $\kappa$. Using this, we can split the central-path equations into the residuals
\begin{equation}\label{eq:jn_splitting}
    G_\mu(u,s)\coloneq
    \begin{bmatrix}
        Q(u)+\mu u-(0,s,\mu/\tau)\\
        0
    \end{bmatrix},
    \qquad
    F_\mu(u,s)\coloneq
    \begin{bmatrix}
        0\\
        s\circ y-\mu e
    \end{bmatrix},
\end{equation}
so that \eqref{eq:central_path} is equivalent to
$F_\mu(u,s)+G_\mu(u,s)=0$. Linearizing $F_\mu$ at $(u_k,s_k)$ while
retaining $G_\mu$ gives
\begin{subequations}\label{eq:Josephy-Newton}
    \begin{align}
        Q(u_{k+1})+\mu u_{k+1}
        &=(0,s_{k+1},\mu/\tau_{k+1}),
        \label{eq:G_wk+1}\\
        s_{k+1}\circ y_k+s_k\circ y_{k+1}
        &=\mu e+s_k\circ y_k.
        \label{eq:jordan_linearization}
    \end{align}
\end{subequations}
Let $D_k\in\R^{m\times m}$ be diagonal, with
\begin{equation}\label{eq:D_k}
    (D_k)_{ii}\coloneq
    \begin{cases}
        0, & i=1,\ldots,z,\\[4pt]
        \dfrac{(s_k)_i}{(y_k)_i}, & i=z+1,\ldots,m.
    \end{cases}
\end{equation}
Note that $D_k\succeq0$ and consequently the block map $u\mapsto(0,D_k y,0)$ is monotone.
Since $D_k y_k=s_k$, equation~\eqref{eq:jordan_linearization}
gives the slack-recovery formula
\begin{equation}\label{eq:slack_recovery}
    s_{k+1}
    =\mu y_k^{-1}-D_k(y_{k+1}-y_k).
\end{equation}
Substituting into \eqref{eq:G_wk+1} eliminates $s_{k+1}$,
leaving
\begin{equation}\label{eq:extended_cp}
        Q(u_{k+1})+\mu u_{k+1}+(0,D_k y_{k+1},0)
       =(0,s_k+\mu y_k^{-1},\mu/\tau_{k+1}),
\end{equation}
with $\kappa_{k+1}=\mu/\tau_{k+1}$.

To accommodate the predictor--corrector scheme below, we
consider the more general system
\begin{equation}\label{eq:generic_system}
        Q(u_{k+1})+\rho u_{k+1}+(0,D_k y_{k+1},0)
        =r+(0,0,\mu/\tau_{k+1}),
\end{equation}
over variables $u_{k+1} = (x_{k+1}, y_{k+1}, \tau_{k+1})$,
where the shift $\rho>0$ may differ from the target $\mu\geq0$,
and $r\coloneq(r_x,r_y,r_\tau)\in\R^{n+m+1}$.
For $\mu>0$, equation~\eqref{eq:extended_cp} is the instance
$\rho=\mu$ and $r=(0,s_k+\mu y_k^{-1},0)$.
The following result gives a closed-form solution to \eqref{eq:generic_system}, adapted from \cite{odonoghue2021}.

\begin{proposition}[Closed-form solution to Eq.~\eqref{eq:generic_system}]
    \label{thm:resolvent}
    Fix a shift $\rho>0$, a target $\mu\geq0$, and a
    right-hand side $r$, and assume $D_k\succeq0$.
    We omit iteration subscripts on the temporary solution vectors.
    Let $x'$, $x''$, $y'$, and
    $y''$ solve the following symmetric quasi-definite system
    with two right-hand sides
    \begin{equation}\label{eq:qd_system}
        \begin{bmatrix}
            P+\rho I & A^\top\\
            A & -D_k-\rho I
        \end{bmatrix}
        \begin{bmatrix}
            x' & x''\\
            y' & y''
        \end{bmatrix}
        =
        \begin{bmatrix}
            r_x & -c\\
            -r_y & b
        \end{bmatrix},
    \end{equation}
    and define
    \begin{align*}
        \theta_2
        &\coloneq
        \rho-c^\top x''-b^\top y''
        -(x'')^\top P x'',\\
        \theta_1
        &\coloneq
        -c^\top x'-b^\top y'
        -2(x')^\top P x''-r_\tau,\\
        \theta_0
        &\coloneq
        (x')^\top P x'+\mu.
    \end{align*}
    Then
    \[
    \theta_2
    =
    \rho\bigl(
        1+\|x''\|^2+\|y''\|^2
    \bigr)
    +(y'')^\top D_k y''
    \geq\rho>0,
    \qquad
    \theta_0\geq\mu\geq0.
    \]
    Equation~\eqref{eq:generic_system} admits a solution with
    $\tau_{k+1}>0$ if and only if $\theta_0>0$ or $\theta_1<0$.
    If $\theta_0=0$ and $\theta_1\geq0$,
    we necessarily have $\mu=0$, and the formula returns
    $\tau_{k+1}=0$. This lies outside the domain of $Q$ but
    is a solution for the extended operator $\mathcal{Q}$.
    In either case the solution is given by
    \[
        \tau_{k+1}
        =
        \frac{-\theta_1+
        \sqrt{\theta_1^2+4\theta_2\theta_0}}{2\theta_2},
    \]
    and
    \[
        x_{k+1}=x'+\tau_{k+1}x'',
        \qquad
        y_{k+1}=y'+\tau_{k+1}y''.
    \]
\end{proposition}
The proof is given in Appendix~\ref{app:ipm}. 
Solving \eqref{eq:qd_system} requires a single $LDL^\top$
factorization and two back-solves. Since the coefficient
matrix is symmetric quasi-definite, an $LDL^\top$
factorization exists under any symmetric permutation
\cite{vanderbei1995}. In a practical solver $D_k$ and $\rho$
changes from one iteration to the next; however, the sparsity pattern is fixed throughout. This permits a fill-reducing ordering
chosen once at the start of the solve and pivot-free
factorization (in exact arithmetic) thereafter. Finite-precision
safeguards are described in Section~\ref{sec:implementation}.

\subsection{Mehrotra-style predictor--corrector}
\label{subsec:mehrotra}

To follow the central path, we adapt the Mehrotra
predictor--corrector algorithm~\cite{mehrotra1992implementation},
widely used for symmetric cones
\cite{goulart2024,vandenberghe2010cvxopt,domahidi2013ecos}.
In these schemes each iteration first computes a \emph{predictor} trial point
targeting zero complementarity, $\mu=0$.
The maximum step length towards this point that preserves cone
membership determines the barrier target for the \emph{corrector},
which balances centering against complementarity reduction.
The iterate is then updated using a damped corrector step.

In standard formulations, changing the barrier target at a fixed
iterate affects only the right-hand side, allowing the predictor
and corrector to share a matrix factorization.
For our central path, the barrier parameter also determines the
diagonal regularization in the coefficient matrix.
We therefore derive a surrogate procedure that holds this shift
fixed across the two solves, preserving factorization reuse
while retaining the regularization provided by the barrier parameter.

For the step construction, let $(u_k,s_k)$ denote the current
iterate after the normalization described below, and define its
average complementarity by
\[
    \mu_k\coloneq\frac{s_k^\top y_k}{\nu},
\]
where we assume throughout this section that $\nu=m-z>0$.
\subsubsection{Surrogate formulation}

To preserve factorization reuse, we use a common shift $\mu_k > 0$
for both predictor and corrector steps. Let $\mu\in[0,\mu_k]$ denote the target barrier parameter for the step. Starting from \eqref{eq:extended_cp} with this target, we introduce a fixed slack-block offset $\delta$ for the second-order correction. Adding $(\mu_k-\mu)u_{k+1}$ to both sides then gives
\begin{equation}\label{eq:split_system}
    Q(u_{k+1})+\mu_k u_{k+1}+(0,D_k y_{k+1},0)
    =(\mu_k-\mu)u_{k+1}
    +(0,s_k+\mu y_k^{-1}+\delta,\mu/\tau_{k+1}).
\end{equation}
Here $\delta$ is zero on equality rows, and $\delta=0$ recovers
the original subproblem~\eqref{eq:extended_cp}.
The term $(\mu_k-\mu)u_{k+1}$ on the right-hand side still depends
on the unknown trial point. We replace this term by
$(\mu_k-\mu)\tilde u$, where $\tilde u$ is a fixed surrogate point,
giving
\begin{equation}\label{eq:surrogate_system}
        Q(u_{k+1})+\mu_k u_{k+1}+(0,D_k y_{k+1},0)
        =(\mu_k-\mu)\tilde u +(0,s_k+\mu y_k^{-1}+\delta,\mu/\tau_{k+1}).
\end{equation}
This is an instance of \eqref{eq:generic_system} with
$\rho=\mu_k$ and
\begin{equation}\label{eq:surrogate_rhs}
    r=(\mu_k-\mu)\tilde u
      +(0,s_k+\mu y_k^{-1}+\delta,0).
\end{equation}
The surrogate system is a single proximal-point step for
the fixed-target subproblem \eqref{eq:split_system}.
For $0<\mu\leq\mu_k$, its solution satisfies
\[
    \|u_{k+1}-u^\star\|
    \leq
    \left(1-\frac{\mu}{\mu_k}\right)
    \|\tilde u-u^\star\|,
\]
where $u^\star$ solves~\eqref{eq:split_system} and $u_{k+1}$
solves~\eqref{eq:surrogate_system} with the same offset $\delta$;
see Proposition~\ref{prop:contraction} in Appendix~\ref{app:ipm}.
This bound motivates choosing a surrogate anchor close to $u^\star$.
We use the current iterate as the surrogate point for the predictor and the predictor trial point as the surrogate point for the corrector.

\subsubsection{Shared linear algebra}

Under the surrogate scheme described above the coefficient matrix depends only on $(s_k,y_k,\mu_k)$ for fixed problem data, so a single sparse factorization of
\begin{equation}\label{eq:matrix_to_factorise}
    K_k \coloneq \begin{bmatrix}
        P+\mu_k I & A^\top\\
        A & -D_k-\mu_k I
    \end{bmatrix}
\end{equation}
serves both step computations.
Moreover, the second right-hand side,
$(-c,b)$, is independent of the target $\mu$, surrogate
$\tilde u$ and offset $\delta$. After factorization, we
therefore solve
\begin{equation}\label{eq:cached_data_solve}
    K_k
    \begin{bmatrix}
        x''\\
         y''
    \end{bmatrix} =
    \begin{bmatrix}
         -c\\
         b
    \end{bmatrix}
\end{equation}
once and cache $(x'',y'')$ for both the predictor
and corrector. Each step then requires only one new solve
with its own right-hand side $(r_x,-r_y)$ to obtain
$(x',y')$. These vectors and the cached pair
determine the scalar quadratic and the trial point through
the closed-form step in Proposition~\ref{thm:resolvent}.

Thus, each iteration uses
one sparse factorization of \eqref{eq:matrix_to_factorise} and three solves with the factors: one for the shared right-hand side and one for each step,
not counting iterative refinement. Although $(-c,b)$ is fixed throughout the solve, its solution must be recomputed whenever the coefficient
matrix changes, in particular at each new outer iteration.

\subsubsection{Normalization and parameter selection}

At the start of each iteration, we rescale $u_k$ and $s_k$
by a common positive factor so that
$\|u_k\|=\sqrt{\nu+1}$, matching the central-path scale from
Proposition~\ref{prop:sphere}. Because the underlying embedding
is homogeneous, shrinking the iterates reduces $s_k^\top y_k$
without changing the recovered primal--dual point
$(x_k/\tau_k,y_k/\tau_k,s_k/\tau_k)$. Normalization prevents
this change of scale from being mistaken for progress by
the predictor--corrector heuristic. We apply the rescaling
\begin{equation}\label{eq:iterate_normalisation}
    (u_k,s_k)\gets
    \frac{\sqrt{\nu+1}}{\|u_k\|}(u_k,s_k).
\end{equation}
This common scaling also preserves cone membership.
After rescaling, we evaluate $\mu_k$ as above and compute a
\emph{predictor} direction
$(\Delta u_k^{\mathrm p},\Delta s_k^{\mathrm p})$ using the zero target $\mu=0$.
For a given direction, define the maximum step length preserving
cone membership
\begin{equation}\label{eq:mehrotra_step_size}
    \alpha(\Delta s,\Delta y)
    \coloneq
    \max\left\{
        \alpha'\in[0,1]:
        (s_k+\alpha'\Delta s,y_k+\alpha'\Delta y)
        \in\mathcal K\times\mathcal K^*
    \right\}.
\end{equation}
We write
$\alpha_k^{\mathrm p}
\coloneq \alpha(\Delta s_k^{\mathrm p},\Delta y_k^{\mathrm p})$
and
$\alpha_k^{\mathrm c}
\coloneq \alpha(\Delta s_k^{\mathrm c},\Delta y_k^{\mathrm c})$
for the predictor and corrector step lengths.
We use $\alpha_k^{\mathrm p}$ without additional damping to evaluate
the centering heuristic. The accepted corrector step is multiplied
by a damping factor
$\gamma\in(0,1)$ to preserve strict interiority.
No additional restriction involving $\Delta\tau$ is needed:
the trial value of $\tau$ is nonnegative
(Proposition~\ref{thm:resolvent}), so damping also preserves $\tau>0$.
Following \cite[Section~3.2]{goulart2024}, we choose
the centering parameter using the cubic step-length heuristic
\begin{equation}\label{eq:centering}
    \sigma_k \coloneq (1 - \alpha_k^{\mathrm p})^3.
\end{equation}
Since $\alpha_k^{\mathrm p}\in[0,1]$, we have $\sigma_k\in[0,1]$.
A predictor step length close to one gives $\sigma_k$ close to zero,
whereas a short feasible step favors stronger centering.
The \emph{corrector} then targets $\mu = \sigma_k\mu_k$ and
incorporates a second-order complementarity correction.
The full predictor trial and its directions are retained for the
corrector computation.

\subsubsection{Predictor, corrector, and update}

For the predictor, we choose the current iterate as the
surrogate, $\tilde u=u_k$, and solve
\eqref{eq:surrogate_system} with $\mu=0$ and $\delta=0$
to obtain $u_{k+1}^{\mathrm p}$. The slack is recovered
from \eqref{eq:slack_recovery}.
For the corrector, Mehrotra's second-order term modifies
the slack-recovery equation to
\begin{equation}\label{eq:corrector_slack}
    s_{k+1}^{\mathrm c}
    =
    \mu y_k^{-1}
    -D_k(y_{k+1}^{\mathrm c}-y_k)
    -y_k^{-1}\circ
        (\Delta s_k^{\mathrm p}\circ\Delta y_k^{\mathrm p}).
\end{equation}
We therefore solve \eqref{eq:surrogate_system} with
\[
    \tilde u=u_{k+1}^{\mathrm p},\qquad
    \mu=\sigma_k\mu_k,\qquad
    \delta=-y_k^{-1}\circ
        (\Delta s_k^{\mathrm p}\circ\Delta y_k^{\mathrm p}),
\]
using the full predictor trial, without damping or renormalization,
as the anchor. We recover $s_{k+1}^{\mathrm c}$ from
\eqref{eq:corrector_slack}.
After computing $\alpha_k^{\mathrm c}$ from
\eqref{eq:mehrotra_step_size}, we update
\[
    (u_{k+1},s_{k+1})
    =
    (u_k,s_k)+\gamma\alpha_k^{\mathrm c}
    (\Delta u_k^{\mathrm c},\Delta s_k^{\mathrm c}).
\]
The new iterate is normalized at the start of the next iteration.

To see the role of the second-order correction, expand the
complementarity products at the full corrector trial.
Equation~\eqref{eq:corrector_slack} gives
\begin{equation}\label{eq:corrector_complementarity_error}
    s_{k+1}^{\mathrm c}\circ y_{k+1}^{\mathrm c}
    -\sigma_k\mu_k e
    =
    \Delta s_k^{\mathrm c}\circ\Delta y_k^{\mathrm c}
    -
    \Delta s_k^{\mathrm p}\circ\Delta y_k^{\mathrm p}.
\end{equation}
Thus, the predictor cross-product estimates the quadratic
contribution omitted by linearizing complementarity, and
the corrector leaves only the mismatch between that estimate
and its own quadratic contribution. This can improve the
accuracy of large trial steps.

\begin{algorithm}
    \caption{Modified Mehrotra predictor--corrector}
    \label{alg:ipm}
    \begin{algorithmic}[1]
        \Require Problem \eqref{eq:primal_qp} with $\nu>0$;
        initial point
        $(u_0,s_0)\in
        \operatorname{int}(\mathcal C)\times
        \operatorname{relint}(\mathcal K)$;
        damping factor $\gamma\in(0,1)$.

        \For{$k=0,1,2,\ldots$}
            \State Rescale $(u_k,s_k)\gets \frac{\sqrt{\nu+1}}{\|u_k\|}(u_k,s_k)$.
            \State Stop if the termination criteria of
            Section~\ref{subsec:termination} are satisfied.
            \State $\mu_k\gets s_k^\top y_k/\nu$;
            form $D_k$ using \eqref{eq:D_k}.

            \Statex \smallskip\textbf{Shared linear algebra}
            \State Form and factorize
            $K_k\gets
            \begin{bmatrix}
                P+\mu_k I & A^\top\\
                A & -D_k-\mu_k I
            \end{bmatrix}$.
            \State Solve
            $K_k\binom{x''}{y''}=\binom{-c}{b}$
            and cache $(x'',y'')$.

            \Statex \smallskip\textbf{Predictor}
            \State $\mu\gets0$,
            $r\gets\mu_k u_k+(0,s_k,0)$.
            \State Solve
            $K_k\binom{x'}{y'}=\binom{r_x}{-r_y}$.
            \State Compute $\tau_{k+1}^{\mathrm p}$ from the
            scalar root formula in Proposition~\ref{thm:resolvent}.
            \State $u_{k+1}^{\mathrm p}\gets
            (x'+\tau_{k+1}^{\mathrm p}x'',
             y'+\tau_{k+1}^{\mathrm p}y'',
             \tau_{k+1}^{\mathrm p})$.
            \State $s_{k+1}^{\mathrm p}\gets
            -D_k(y_{k+1}^{\mathrm p}-y_k)$.
            \State
            $\Delta u_k^{\mathrm p}\gets
            u_{k+1}^{\mathrm p}-u_k$,
            $\Delta s_k^{\mathrm p}\gets
            s_{k+1}^{\mathrm p}-s_k$.
            \State Compute $\alpha_k^{\mathrm p}$
            using \eqref{eq:mehrotra_step_size}.
            \State $\sigma_k\gets(1-\alpha_k^{\mathrm p})^3$.

            \Statex \smallskip\textbf{Corrector}
            \State $\mu\gets\sigma_k\mu_k$,
            $\delta\gets-y_k^{-1}\circ
            (\Delta s_k^{\mathrm p}\circ\Delta y_k^{\mathrm p})$.
            \State $r\gets(\mu_k-\mu)u_{k+1}^{\mathrm p}
            +(0,s_k+\mu y_k^{-1}+\delta,0)$.
            \State Solve
            $K_k\binom{x'}{y'}=\binom{r_x}{-r_y}$.
            \State Compute $\tau_{k+1}^{\mathrm c}$ from the
            scalar root formula in Proposition~\ref{thm:resolvent}.
            \State $u_{k+1}^{\mathrm c}\gets
            (x'+\tau_{k+1}^{\mathrm c}x'',
             y'+\tau_{k+1}^{\mathrm c}y'',
             \tau_{k+1}^{\mathrm c})$.
            \State $s_{k+1}^{\mathrm c}\gets
            \mu y_k^{-1}-D_k(y_{k+1}^{\mathrm c}-y_k)+\delta$.
            \State
            $\Delta u_k^{\mathrm c}\gets
            u_{k+1}^{\mathrm c}-u_k$,
            $\Delta s_k^{\mathrm c}\gets
            s_{k+1}^{\mathrm c}-s_k$.

            \Statex \smallskip\textbf{Update}
            \State Compute $\alpha_k^{\mathrm c}$
            using \eqref{eq:mehrotra_step_size}.
            \State $(u_{k+1},s_{k+1})\gets
            (u_k,s_k)+\gamma\alpha_k^{\mathrm c}
            (\Delta u_k^{\mathrm c},\Delta s_k^{\mathrm c})$.
        \EndFor
    \end{algorithmic}
\end{algorithm}

\paragraph{Proximal interpretation.}
The surrogate formulation \eqref{eq:surrogate_system} can be interpreted as dynamic
proximal regularization towards a step-dependent anchor.
For a target $\mu\leq\mu_k$, define the residual function of \eqref{eq:split_system} as
$$S_k(u) \coloneq Q(u) + \mu u +(0, D_ky,0) - (0,s_k+\mu y_k^{-1} + \delta, \mu / \tau).$$
Rearranging \eqref{eq:split_system} and \eqref{eq:surrogate_system} gives
$$S_k(u^\star) = 0, \quad S_k(u_{k+1}) + (\mu_k - \mu)(u_{k+1} - \tilde u) = 0,$$
where the latter is equivalent to 
$$u_{k+1} = (I + (\mu_k - \mu)^{-1}S_k)^{-1}(\tilde u).$$
The surrogate system is therefore equivalent to one iteration of the proximal point method for solving the monotone equation $S_k(u^\star) = 0$, with \textit{proximal weight} $(\mu_k - \mu)^{-1}$ and \textit{anchor} $\tilde u$. The distance between $u_{k+1}$ and the anchor is controlled by the weight, leading to regularization in the step directions. In the predictor step, $\mu=0$ and $\tilde u=u_k$, so
the proximal weight discourages
large departures from the current iterate.
In the corrector step, the anchor is the predictor trial
point $u_{k+1}^{\mathrm p}$ and the regularization weight
is $(1-\sigma_k)\mu_k$.
When the feasible predictor step is short, $\sigma_k\approx1$ and the corrector
emphasizes centering at approximately the current barrier
parameter, with low regularization.
Conversely, when the predictor step is close to one,
$\sigma_k\approx0$ and the corrector targets a more aggressive
reduction in the barrier parameter, while also retaining
high regularization towards the predictor trial point.
In both steps, the proximal regularization weight vanishes
as $\mu_k\to0$.

\paragraph{Algorithm summary.}
Algorithm~\ref{alg:ipm} makes the three linear solves explicit.
In each step, $r=(r_x,r_y,r_\tau)$ is the right-hand side
of \eqref{eq:generic_system}. The scalar root formula in
Proposition~\ref{thm:resolvent} uses the fixed shift
$\rho=\mu_k$, the current target $\mu$, and the scalar
component $r_\tau$.
Both steps reuse the factorization of $K_k$ and the cached
$(x'',y'')$.

\section{Implementation}\label{sec:implementation}

Our software package \cpot implements the interior-point method developed in this paper
for convex quadratic and linear programs, with
$\mathcal{K}=\{0\}^{z}\times\R_{\geq0}^{m-z}$.
The solver is written in Rust, provides a Python interface, and is freely
available at
\begin{equation}\label{eq:source_code}
    \text{\url{https://github.com/cpo3-org/cpo3}}.
\end{equation}
The name \cpot{} is a play on `conic programming in Rust',
with `O3' alluding to iron oxide ($\mathrm{Fe_2O_3}$). The following subsections describe the initialization, numerical linear algebra, and the safeguards used to implement Algorithm~\ref{alg:ipm}. By default, we use damping factor $\gamma = 0.999$.

\subsection{Presolve and data equilibration}

Before constructing the embedding, we optionally apply presolve to reduce the problem
dimensions by eliminating variables and redundant constraints. We implement
standard presolve methods described in~\cite{cederberg2026presolving,gleixner2023papilo,achterberg2020presolve},
with modifications where necessary to accommodate quadratic objectives.
Although presolve is disabled in our main numerical experiments, we
summarize its implementation in Appendix~\ref{app:presolve} for completeness.

We scale the problem data using symmetric Ruiz equilibration~\cite{ruiz2001},
a technique commonly used in modern
solvers~\cite{stellato2020,garstka2021,goulart2024,schwan2023piqp,ocpb:16,applegate2021practical,applegate2025pdlp}.
This scaling is applied whether presolve is enabled or not. Specifically, we equilibrate
the symmetric matrix
\[
 B \coloneq \begin{bmatrix}
    P & A^\top & c \\
    A & 0 & b \\
    c^\top & b^\top & 0
 \end{bmatrix}.
\]
Starting with $D=I$, we compute the row infinity norms
$r_i=\max_j|(DBD)_{ij}|$ and update the diagonal scaling factors simultaneously by
$D_{ii}\leftarrow D_{ii}/\sqrt{r_i}$ when $r_i>0$.
By default, we perform at most two equilibration sweeps,
stopping early if all nonzero row norms are within $10^{-2}$
of one. The interior-point method operates on the equilibrated problem data, and recovered solutions are mapped back to the original
coordinates.

\subsection{Default initialization}\label{subsec:initialization}

By default, we initialize at the large-parameter endpoint $u_\infty$
from Proposition~\ref{prop:large_mu_endpoint}. Specifically, we set
\[
    x_0=0,\qquad y_0=s_0=e,\qquad \tau_0=1,
\]
where $e$ equals one on inequality rows and zero on equality rows.
This initialization requires no linear solve and satisfies
$\|u_0\|^2=\nu+1$, with initial complementarity
$s_0^\top y_0/\nu=1$ when $\nu>0$.

\subsection{Warm-start initialization}\label{subsec:warm-start}

A primal--dual solution of a related problem provides an alternative to the
default initialization. After any presolve reductions and data equilibration,
we construct three strictly interior candidates by
flooring the inequality slacks and multipliers at $10^{-12}$, $10^{-6}$,
and $10^{-3}$, respectively. For each candidate, we set $\tau=1$ and rescale
all homogeneous coordinates so that $\|(x,y,\tau)\|=\sqrt{\nu+1}$.
We compare the candidates using the proximity measure $\eta_\mu(u)$ defined
in~\eqref{eq:eta_def}, evaluated at each candidate's complementarity
$\mu=s^\top y/\nu$, and select the candidate with the smallest value.
This selection rule is motivated by the local-norm error bound in
Theorem~\ref{thm:barrier_error_bound}, although we do not enforce $\eta_\mu(u)<1$ for any candidate warm start.

Warm-starting interior-point methods is known to be difficult: the optimal
solution of a nearby problem lies close to the boundary of the cone and is
poorly centered, so the iterates can make slow progress at
first~\cite{gondzio2008unblocking}.
There is a highly active line of research investigating effective warm-starting for
IPMs \cite{yildirim2002warm,gondzio2002reoptimization,john2008implementation,benson2007exact,skajaa2013warmstarting,raghunathan2022homogeneous,chen2025warmstarting}.
For homogeneous paths anchored at the starting point, a warm start also
determines where the path meets the solution ray
(Figure~\ref{fig:central_paths}). Our central path does not
depend on the starting point: a warm start cannot change the path, only the
initial distance to it, which $\eta_\mu(u)$ measures
(Theorem~\ref{thm:barrier_error_bound}). Flooring the slacks and multipliers
plays the role of moving the previous optimum off the boundary, and
$\eta_\mu$ selects among the floors.

\subsection{KKT factorization}\label{subsec:kkt_factorization}

Rather than computing a factorization of the matrix in \eqref{eq:matrix_to_factorise}, we first eliminate variables associated
with inequality rows of $A$ containing exactly one nonzero entry,
such as simple bound constraints. After symmetrically permuting $M_k$
to place these variables last, a system with right-hand side $(r_R,r_E)$
has the form
\[
    \underbrace{\begin{bmatrix}
        B_k & C \\
        C^\top & -\Theta_k
    \end{bmatrix}}_{\text{permuted} \ M_k}
    \begin{bmatrix} z_R \\ z_E \end{bmatrix}
    = \begin{bmatrix} r_R \\ r_E \end{bmatrix},
\]
where $R$ and $E$ denote retained and eliminated variables, respectively.
Each column of $C$ corresponds to an eliminated row of $A$ and contains
its single nonzero entry, while $\Theta_k$ is positive diagonal.
Taking the Schur complement gives the reduced system
\[
    \left(B_k+C\Theta_k^{-1}C^\top\right)z_R
    =r_R+C\Theta_k^{-1}r_E.
\]
The contribution $C\Theta_k^{-1}C^\top$ is diagonal, so elimination
introduces no fill-in. The eliminated variables are recovered by
\[
    z_E=\Theta_k^{-1}\left(C^\top z_R-r_E\right).
\]
This elimination is enabled by default, independently of presolve.

\subsection{Iterative refinement}\label{subsec:refinement}

Both predictor and corrector use $\rho=\mu_k$ in~\eqref{eq:qd_system}. Denote by
$K_k$ the matrix in \eqref{eq:matrix_to_factorise} to be factorized.
Very close to the optimum, the small regularization supplied by
$\mu_k$ can make accurate factorization difficult
\cite{gill1996stability}.
In those cases, we factor an approximation $M_k\approx K_k$ with adaptive diagonal regularization (explained in the sequel). We apply corrections
only if the residual test below fails, using
\[
    z^{j+1} \gets z^j + M_k^{-1}r^j,
    \qquad r^j \coloneq r-K_kz^j.
\]
This is preconditioned Richardson iteration with unit step
length~\cite{golub2013matrix,saad2003iterative}.
For $r=(-c,b)$, the seed is the previous outer iteration's solution
when available; otherwise, and for predictor and corrector solves,
it is zero. One application of the update produces $z^0$ and counts
as the initial solve; subsequent corrections count as refinement steps.
We stop when
\[
    \max_{1\leq i\leq n+m}
    \frac{|r^j_i|}
         {\bigl(|r|+|K_k||z^0|\bigr)_i}
    \leq 10^{-7},
\]
with absolute values taken componentwise. We allow at most $30$
refinement steps, stopping early on divergence or when the observed
convergence rate predicts failure to reach tolerance within the
remaining budget. If the tolerance is not reached, we retain the
iterate with the smallest scaled residual.

\subsection{Adaptive diagonal regularization}\label{subsec:regularization}

We construct $M_k$ by imposing a diagonal floor $\zeta_k>0$ on
$K_k$, preserving signs and leaving off-diagonal entries unchanged:
\[
\begin{aligned}
(M_k)_{ii}
&=\max\{P_{ii}+\mu_k,\zeta_k\},
&&i=1,\ldots,n,\\
(M_k)_{n+i,n+i}
&=-\max\{(D_k)_{ii}+\mu_k,\zeta_k\},
&&i=1,\ldots,m.
\end{aligned}
\]
A larger floor can improve factorization stability but may slow
refinement by increasing the discrepancy $M_k-K_k$.
We adapt $\zeta_k$ using the residual after each computed correction
$\Delta z^j\coloneq z^{j+1}-z^j$, which decomposes as
\[
    r^{j+1}
    =
    \underbrace{(M_k-K_k)\Delta z^j}_{\text{regularization residual}}
    +
    \underbrace{r^j-M_k\Delta z^j}_{\text{solve residual}}.
\]
When refinement stops above tolerance, we compare these contributions
for the last correction. For each contribution, we divide its
componentwise magnitude by $|r|+|K_k||z^0|$ and record the maximum.
After the three linear solves, we update the floor for the next
outer iteration. We increase it if the solve residual dominates
in at least one unsuccessful solve and the regularization residual
dominates in none. We decrease it under the converse condition,
or if all solves meet tolerance; otherwise, we leave it unchanged.

\subsection{Linearization fallback for the homogeneous step}

Although Proposition~\ref{thm:resolvent} guarantees a finite,
nonnegative $\tau$ in exact arithmetic, its numerical evaluation
can fail. We then use the following fallbacks, reusing the
KKT solutions from~\eqref{eq:qd_system} to define
\[
    u(\tau)\coloneq
    (x'+\tau x'',\,y'+\tau y'',\,\tau).
\]
For the current predictor or corrector step, write the
$\tau$-row of~\eqref{eq:generic_system}, multiplied by $\tau$, as
\[
    g(u)\coloneq
    \mu_k\tau^2-\tau(c^\top x+b^\top y+r_\tau)
    -x^\top Px-\mu=0,
\]
where $\mu$ and $r_\tau$ are the target parameter and scalar
right-hand side for that step.
We first replace $x^\top Px$ by its linearization
$2x_k^\top Px-x_k^\top Px_k$ and solve the resulting scalar
quadratic in $\tau$. This removes the quadratic dependence
of its coefficients on the approximate KKT solutions.
If this quadratic also fails numerically, we instead linearize
$g$ at the current iterate $u_k=(x_k,y_k,\tau_k)$:
\[
    g(u_k)+\nabla g(u_k)^\top\bigl(u(\tau)-u_k\bigr)=0.
\]
This equation is linear in $\tau$. If it is numerically singular
or yields a nonfinite value, we set $\tau=\tau_k$.
For either fallback, we clip $\tau$ to
$[0.1\tau_k,10\tau_k]$ before recovering
$x=x'+\tau x''$ and $y=y'+\tau y''$.
The original quadratic solution is not clipped.

\subsection{Termination and solution recovery}\label{subsec:termination}

\paragraph{Optimality.}
For $\tau>0$, we form a candidate solution in the original problem
coordinates by undoing equilibration, dividing by $\tau$, and reversing
any presolve reductions. Below, $(x,y,s)$ denotes this candidate, with
primal and dual objective values
\[
p \coloneq \frac{1}{2} x^\top Px+c^\top x \quad \text{and} \quad
d \coloneq -\frac{1}{2} x^\top Px-b^\top y,
\]
respectively (using $P\xi = Px$). We declare optimality when both feasibility tests and
either gap test are satisfied. The feasibility tests are
\begin{subequations}\label{eq:feasibility_tests}
\begin{align}
 \frac{\|Ax+s-b\|_\infty}
 {1+\max\{\|Ax\|_\infty,\|s\|_\infty,\|b\|_\infty\}}
 &<\varepsilon_{\mathrm{feas}},\label{eq:primal_feasibility}\\
 \frac{\|Px+A^\top y+c\|_\infty}
 {1+\max\{\|Px\|_\infty,\|A^\top y\|_\infty,\|c\|_\infty\}}
 &<\varepsilon_{\mathrm{feas}},\label{eq:dual_feasibility}
\end{align}
\end{subequations}
and the gap tests are
\[
 |p-d|<\varepsilon_{\mathrm{abs}}, \qquad
 \frac{|p-d|}{\max\{1,\min(|p|,|d|)\}}<\varepsilon_{\mathrm{rel}}.
\]
All three tolerances default to $10^{-8}$.

\paragraph{Infeasibility.}
A large ratio $\kappa/\tau$ indicates that the iterate may be approaching
an infeasibility certificate. For $\nu>0$, we test for such certificates
after a step only when
\[
 \frac{\kappa}{\tau}=\frac{s^\top y}{\nu\tau^2}>10^9.
\]
We undo equilibration on the homogeneous ray without dividing by $\tau$.
Using $(x,y,s)$ for this ray and the data of the problem entering
equilibration, primal infeasibility requires
\[
 b^\top y<-\varepsilon_{\mathrm{inf,abs}},\qquad
 \frac{\|A^\top y\|_2}
 {\max\{|b^\top y|,\|y\|_2,10^{-15}\}}
 <\varepsilon_{\mathrm{inf,rel}}.
\]
Dual infeasibility requires
\begin{gather*}
 c^\top x<-\varepsilon_{\mathrm{inf,abs}},\\
 \frac{\|Ax+s\|_2}
 {\max\{|c^\top x|,\|x\|_2+\|s\|_2,10^{-15}\}}
 <\varepsilon_{\mathrm{inf,rel}},\\
 \frac{\|Px\|_2}
 {\max\{|c^\top x|,\|x\|_2,10^{-15}\}}
 <\varepsilon_{\mathrm{inf,rel}}.
\end{gather*}
Both infeasibility tolerances default to $10^{-8}$. We normalize primal and
dual infeasibility certificates by $-b^\top y$ and $-c^\top x$, respectively,
then reverse any presolve transformations to return a certificate for the
original problem.

\section{Numerical experiments}\label{sec:numerical}

Here we evaluate the solver's reliability and runtime on a range of LP and QP benchmarks,
and compare to several other open-source and proprietary solvers.
We also examine its iteration and linear-algebra costs relative to Clarabel, another solver targeting the homogeneous embedding, and finally
assess the benefit of the proposed warm-start heuristic.

\paragraph{Experimental setup.}
We evaluate 675 problems from eight collections:
Maros--M\'esz\'aros (138 instances) \cite{maros1999repository},
QPLIB (19 instances) \cite{furini2019qplib},
MPC (64 instances) \cite{caron2023qpbenchmark},
DC-OPF (68 MATPOWER instances) \cite{zimmerman2011matpower},
NETLIB (93 instances), Kennington (16 instances),
Mittelmann (37 instances) \cite{mittelmann2021}, and
MIPLIB~2017 (240 instances) \cite{gleixner2021miplib}.
We compare \cpot{} against a collection of leading open-source and commercial solvers---Clarabel~0.11.1,
HiGHS~1.15.1, PIQP~0.6.4, Gurobi~13.0.3, MOSEK~11.2.4, and CPLEX~22.2.0.1.
The comparisons against other solvers run on Berkeley's Savio cluster, on nodes with Intel Xeon Gold 6130
(Skylake) processors, 32 cores, and 384 GB of RAM.
Each job reserves one node exclusively and runs up to 16
independent single-threaded solves, each limited to 20 GiB
of memory and 300 seconds of runtime.

\paragraph{Verification and tolerance selection.}
To compare solvers with different stopping criteria, we independently
verify returned candidates using the original problem data.
For $(x,y,s)$, the primal and dual feasibility residuals are the
left-hand sides of~\eqref{eq:primal_feasibility}
and~\eqref{eq:dual_feasibility}. The duality-gap, primal-cone,
and dual-cone residuals are, respectively,
\[
\frac{|p-d|}{1+|p|+|d|},\qquad
\frac{\|s-\Pi_{\mathcal K}(s)\|_\infty}
     {1+\max\{\|Ax\|_\infty,\|s\|_\infty,\|b\|_\infty\}},\qquad
\frac{\|y-\Pi_{\mathcal K^*}(y)\|_\infty}
     {1+\|y\|_\infty},
\]
where $p$ and $d$ are the candidate's primal and dual objective
values, and $\Pi$ denotes Euclidean projection onto the indicated cone.

\begin{figure}
    \centering
    \includegraphics[width=0.9\linewidth]{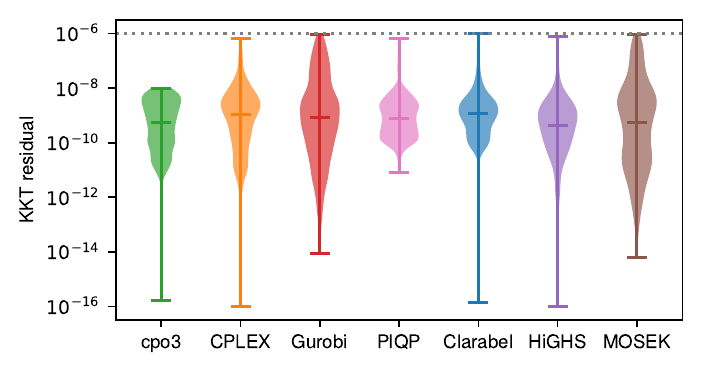}
    \caption{Distribution of the maximum KKT residual for each
    solver's verified trials. The dotted line marks the
    verification threshold $10^{-6}$, and horizontal marks
    within the violins indicate medians. Under these settings the solvers terminate with very similar final accuracies.}
    \label{fig:solver_residuals}
\end{figure}

We set solver tolerances to $10^{-8}$ and accept a candidate as
verified if all five residuals are below $10^{-6}$.
This verification threshold is two orders of magnitude above the
requested tolerance, allowing for differences in residual definitions
and stopping criteria while imposing a common accuracy requirement.
Figure~\ref{fig:solver_residuals} shows, for each solver, the
distribution of the maximum residual per verified trial.
The residuals are generally of similar magnitude across solvers,
indicating comparable solution accuracy among verified trials.

\paragraph{Performance measures.}
For the comparisons in Section~\ref{sec:central-path} and
Appendix~\ref{app:default-settings}, we perform three trials
per solver--problem pair, applying the verification
criterion above to each trial separately. We assign each
trial its measured runtime if verified and infinite runtime
otherwise, then take the median of the three values.
A problem therefore counts as solved if at least two trials are verified.

Tables report verified solve counts and shifted geometric mean runtimes,
\[
  \exp\left(\frac{1}{N}\sum_{i=1}^{N}\log(t_i+10)\right)-10,
\]
where $N$ is the number of problems in the reported group
and $t_i$ is the median runtime in seconds for a solved
problem and 300 seconds otherwise (\emph{i.e.}, the runtime budget).
Times below 0.01 seconds are shown as $<0.01$, and bold entries indicate the best values before rounding.

To complement the aggregate runtime statistics, we use
Dolan--Mor\'e performance profiles~\cite{dolan2002benchmarking}.
For each solver, the profile shows the fraction of problems
whose median runtime is within a given factor of the fastest
median runtime on that problem. At a factor of one, the curve
gives the fraction on which the solver is fastest; its
limiting height is the fraction solved.
Higher curves therefore indicate better performance.

\subsection{Solver performance on feasible problems}\label{sec:central-path}

We use each solver's interior-point algorithm, with presolve and crossover
disabled where available. Overall, \cpot{} solves 631 of 675 problems (93.5\%),
compared with 595 for CPLEX, which has the second-highest solve count. Its
penalized shifted geometric mean runtime is also the lowest: 5.60 seconds,
compared with 8.18 seconds for CPLEX.
\cpot{} also leads in both verified solve count and penalized shifted geometric
mean runtime when LPs and QPs are considered separately, outperforming all tested
open-source and commercial solvers on these measures.

\subsubsection{Quadratic programs}

\begin{figure}[!htbp]
\centering
\includegraphics[width=0.9\linewidth]
{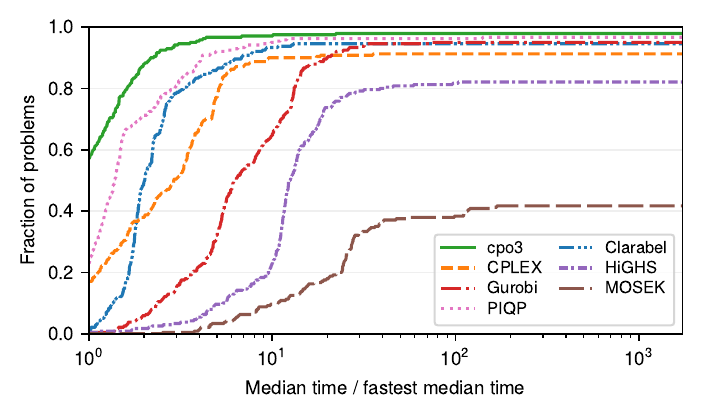}
\caption{Dolan--Mor\'e performance profiles for QPs.}
\label{fig:central-path-qp}
\end{figure}

Table~\ref{tab:central-path-qp} and Figure~\ref{fig:central-path-qp} report the QP results.
\cpot{} solves 235 of 240 QPs, the highest count among the tested solvers,
followed by PIQP with 232. It also achieves the lowest penalized shifted
geometric mean runtime at 1.37 seconds, compared with 2.25 seconds for PIQP.
\cpot{} solves all MPC and DC-OPF instances and is the only solver to solve
every Maros--M\'esz\'aros instance. It also has the lowest penalized shifted
geometric mean runtime on the Maros--M\'esz\'aros and MPC datasets.

MOSEK reports optimality on 116 QPs that fail independent verification, primarily because of excessive stationarity residuals. MOSEK does not target QPs directly, instead it
transforms them into equivalent SOCPs. However, this reformulation can return solutions that have higher than expected residuals.
In diagnostic runs, tightening tolerances recovered verified solutions for some problems but introduced regressions in other problems.

For HiGHS, we use its HiPO interior-point solver~\cite{zanetti2025factorisation}.
Of its 42 QPs without a verified solution, 35 terminate with warnings
indicating stagnation or insufficient accuracy.

\begin{table}[!htbp]
\centering
\small
\setlength{\tabcolsep}{3pt}
\sisetup{detect-weight=true}
\caption{QP results.}
\label{tab:central-path-qp}
\begin{tabular}{@{}l
r S[table-format=3.2]
r S[table-format=3.2]
r S[table-format=<2.2]
r S[table-format=3.2]
@{}}
\toprule
Solver
& \multicolumn{2}{c}{Maros--M\'esz\'aros}
& \multicolumn{2}{c}{QPLIB}
& \multicolumn{2}{c}{MPC}
& \multicolumn{2}{c}{DC-OPF}
\\
& \multicolumn{2}{c}{$N=138$}
& \multicolumn{2}{c}{$N=19$}
& \multicolumn{2}{c}{$N=64$}
& \multicolumn{2}{c}{$N=19$}
\\
\cmidrule(lr){2-3}\cmidrule(lr){4-5}\cmidrule(lr){6-7}\cmidrule(lr){8-9}
 & Solved & {Time (s)} & Solved & {Time (s)} & Solved & {Time (s)} & Solved & {Time (s)} \\
\midrule
\cpot{} & \textbf{138} & \bfseries 0.51 & 14 & 23.18 & \textbf{64} & {\boldmath$<0.01$} & \textbf{19} & 0.61 \\
CPLEX & 133 & 1.69 & 14 & \bfseries 21.56 & \textbf{64} & {$<0.01$} & 8 & 63.02 \\
Gurobi & 130 & 3.44 & \textbf{15} & 29.35 & \textbf{64} & {$<0.01$} & \textbf{19} & 0.48 \\
PIQP & 137 & 1.09 & 12 & 48.28 & \textbf{64} & {$<0.01$} & \textbf{19} & 0.52 \\
Clarabel & 130 & 2.81 & 14 & 31.04 & \textbf{64} & {$<0.01$} & \textbf{19} & \bfseries 0.27 \\
HiGHS & 122 & 5.95 & 11 & 47.36 & 61 & 1.75 & 4 & 143.53 \\
MOSEK & 57 & 65.91 & 1 & 249.77 & 35 & 37.46 & 7 & 77.55 \\
\bottomrule
\end{tabular}
\end{table}

\subsubsection{Linear programs}

\begin{figure}[!htbp]
\centering
\includegraphics[width=0.9\linewidth]
{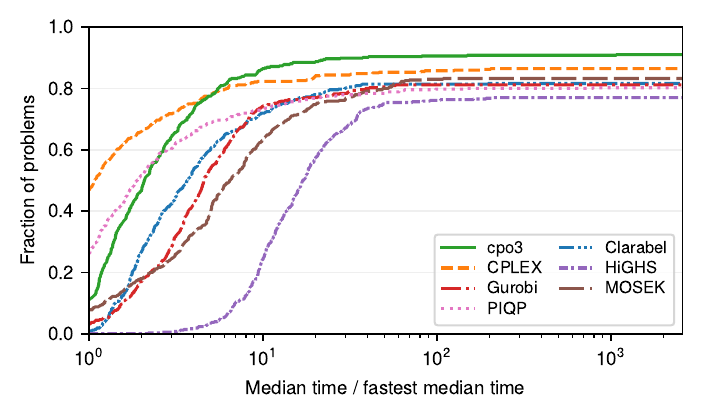}
\caption{Dolan--Mor\'e performance profiles for LPs.}
\label{fig:central-path-lp}
\end{figure}

Table~\ref{tab:central-path-lp} and Figure~\ref{fig:central-path-lp} report the LP results.
\cpot{} solves 396 of 435 LPs, the highest count among the tested solvers,
followed by CPLEX with 376. It also achieves the lowest penalized shifted
geometric mean runtime at 8.58 seconds, compared with 10.98 seconds for CPLEX.

These results show that the method is also very competitive for linear programming.
Beyond solving every DC-OPF instance, it solves all 16 Kennington problems
and all 93 NETLIB problems. On the largest collection, MIPLIB, it verifies
225 of 240 instances, close to CPLEX's leading count of 227 and ahead of
every other tested solver. It also achieves the highest verified solve
count among open-source solvers on Mittelmann.

CPLEX verifies only 21 of the 49 DC-OPF LPs under the benchmark configuration,
despite its strong performance on the other LP collections. On 16 instances,
CPLEX reports that numerical difficulties prevent convergence to the requested
tolerances; these instances fail both the primal-feasibility and duality-gap checks.
The remaining 12 are reported optimal and satisfy the duality-gap criterion,
but eight fail primal feasibility and four fail stationarity. In each case,
our absolute residual for the failed check---the numerator
in~\eqref{eq:primal_feasibility} or~\eqref{eq:dual_feasibility},
respectively---agrees with CPLEX's reported residual within $10^{-10}$.
CPLEX thus accepts these solutions as optimal even though the corresponding
relative residuals exceed our verification tolerance.
With presolve and crossover enabled, however, CPLEX verifies all 49 DC-OPF LPs;
see Appendix~\ref{app:default-settings}.

\begin{table}[!htbp]
\centering
\small
\setlength{\tabcolsep}{3pt}
\sisetup{detect-weight=true}
\caption{LP results.}
\label{tab:central-path-lp}
\begin{adjustbox}{max width=\linewidth}
\begin{tabular}{@{}l
r S[table-format=2.2]
r S[table-format=2.2]
r S[table-format=2.2]
r S[table-format=3.2]
r S[table-format=2.2]
@{}}
\toprule
Solver
& \multicolumn{2}{c}{DC-OPF}
& \multicolumn{2}{c}{NETLIB}
& \multicolumn{2}{c}{Kennington}
& \multicolumn{2}{c}{Mittelmann}
& \multicolumn{2}{c}{MIPLIB}
\\
& \multicolumn{2}{c}{$N=49$}
& \multicolumn{2}{c}{$N=93$}
& \multicolumn{2}{c}{$N=16$}
& \multicolumn{2}{c}{$N=37$}
& \multicolumn{2}{c}{$N=240$}
\\
\cmidrule(lr){2-3}\cmidrule(lr){4-5}\cmidrule(lr){6-7}\cmidrule(lr){8-9}\cmidrule(lr){10-11}
 & Solved & {Time (s)} & Solved & {Time (s)} & Solved & {Time (s)} & Solved & {Time (s)} & Solved & {Time (s)} \\
\midrule
\cpot{} & \textbf{49} & \bfseries 0.10 & \textbf{93} & 0.64 & \textbf{16} & 7.65 & 13 & 149.32 & 225 & 8.81 \\
CPLEX & 21 & 61.64 & 92 & \bfseries 0.46 & \textbf{16} & \bfseries 1.48 & \textbf{20} & \bfseries 112.70 & \textbf{227} & \bfseries 6.95 \\
Gurobi & 42 & 6.44 & 83 & 4.73 & \textbf{16} & 3.83 & 18 & 127.77 & 194 & 16.99 \\
PIQP & \textbf{49} & 0.13 & 92 & 1.25 & 15 & 8.15 & 6 & 220.10 & 188 & 17.57 \\
Clarabel & 44 & 4.28 & 87 & 2.79 & 15 & 13.86 & 9 & 200.51 & 200 & 18.34 \\
HiGHS & 38 & 12.15 & 78 & 8.50 & 12 & 35.10 & 6 & 238.55 & 201 & 24.06 \\
MOSEK & 43 & 5.40 & 85 & 3.71 & 14 & 12.86 & 16 & 117.99 & 204 & 13.03 \\
\bottomrule
\end{tabular}
\end{adjustbox}
\end{table}

\paragraph{Default presolve and crossover settings.}
With the solvers' default presolve and crossover settings enabled where available
(Appendix~\ref{app:default-settings}),
\cpot{} continues to solve more QPs than any of the tested
open-source or commercial solvers, while achieving the lowest penalized shifted
geometric mean runtime: 1.28 seconds, compared with 2.10 seconds for Gurobi, the
next fastest. Presolve also improves \cpot{}'s LP performance, yielding 9
additional verified solves and reducing its LP mean runtime from 8.58 to 6.48
seconds. Across the full corpus, \cpot{} achieves the highest verified solve
count, 640 of 675, and has the lowest penalized shifted geometric mean runtime
among the open-source solvers.

\subsection{Comparison with Clarabel}\label{sec:cpo3-vs-clarabel}

In this section we compare \cpot{} with Clarabel in greater detail, examining solution coverage,
outer iterations, linear-algebra work, and primal-infeasibility detection. Clarabel
is a natural reference solver because both methods operate directly on the
homogeneous embedding. Moreover, both solvers are implemented
in Rust and use the same \verb|faer| linear algebra backend for the LDL factorizations
\cite{Kazdadi2026}.\footnote{By default, Clarabel chooses between QDLDL and
\texttt{faer} for each problem. We retained this behavior in the preceding
experiment; here, we force Clarabel to use \texttt{faer} for all problems to
match \cpot{}.}
We run both solvers with presolve disabled and feasibility and absolute and relative
gap tolerances set to $10^{-8}$. We use the independent verification threshold
$10^{-6}$ described above.

\subsubsection{Solve counts and outer iterations}

Table~\ref{tab:clarabel-coverage} compares verified solve counts on the
675 feasible problems. \cpot{} solves 631 problems and Clarabel solves
581, with 580 solved by both. Of the remaining problems, 51 are solved
only by \cpot{} and one only by Clarabel.
The additional solves for \cpot{} include 27 MIPLIB instances, eight
Maros--M\'esz\'aros instances, and five Mittelmann instances.
Every LP solved by Clarabel is also solved
by \cpot{}.

On the 580 problems verified by both solvers, the median outer iteration
count is 18 for \cpot{} and 14.5 for Clarabel. The difference is concentrated
in LPs, where the medians are 25 and 19; both solvers have a median of 11
iterations on QPs. Nevertheless, \cpot{} requires fewer applications of the
KKT factors, as the reduction in refinement corrections more than offsets
the additional initial solves.

\begin{table}[!htbp]
\centering
\caption{Feasible-problem comparison with Clarabel. Unique solves are relative
to these two solvers. Iteration medians use only the shared verified problems
in each group.}
\label{tab:clarabel-coverage}
\begin{tabular}{@{}llrrrr@{}}
\toprule
Problems & Solver & Solved & Unique & Shared & Median iterations \\
\midrule
All & \cpot{} & 631 & 51 & 580 & 18 \\
 & Clarabel & 581 & 1 & 580 & 14.5 \\
\addlinespace
LP & \cpot{} & 396 & 42 & 354 & 25 \\
 & Clarabel & 354 & 0 & 354 & 19 \\
\addlinespace
QP & \cpot{} & 235 & 9 & 226 & 11 \\
 & Clarabel & 227 & 1 & 226 & 11 \\
\bottomrule
\end{tabular}
\end{table}

\subsubsection{Linear-algebra work}

We count numerical factorizations, KKT solves, and refinement
corrections for each problem, including initialization and unsuccessful
attempts. Each initial solve and refinement correction applies the KKT
factors once, so their combined count gives the total number of factor
applications.

Table~\ref{tab:clarabel-linear-work} reports median counts on the 580
problems verified by both solvers. Although \cpot{} performs more initial
KKT solves at the median, it requires fewer refinement corrections:
one compared with 64 for Clarabel. The median number of factor applications
is consequently lower, at 60 compared with 111.5, with lower medians on
both LPs and QPs. Across all 580 shared problems, \cpot{} uses 36.3\% fewer factor
applications: 67,839 versus 106,575. It performs more numerical
factorizations, 16,787 versus 13,534, consistent with its higher total
outer iteration count.

\begin{figure}[!htbp]
\centering
\includegraphics[width=0.9\linewidth]
{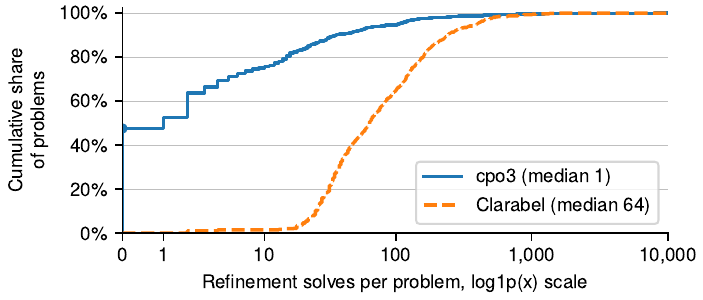}
\caption{Cumulative distributions of total refinement corrections per problem
on the 580 feasible problems verified by both solvers. Each curve
gives the fraction of problems requiring at most the indicated number
of corrections. The horizontal axis uses a $\log(1+x)$ scale to include
zero corrections. Overall \cpot{} requires substantially fewer refinement
corrections than Clarabel.}
\label{fig:clarabel-refinement}
\end{figure}

Figure~\ref{fig:clarabel-refinement} shows how the refinement counts
are distributed across these problems.
\cpot{} requires no refinement corrections on 276 problems (47.6\%),
whereas Clarabel requires at least one on every problem.
At most ten corrections suffice for \cpot{} on 439 problems (75.7\%),
compared with nine (1.6\%) for Clarabel.
Lower refinement is primarily driven by the fact that the matrices
factorized by \cpot{} are naturally regularized by the barrier parameter,
and therefore quasi-definite and stable factorizations always exist (in exact arithmetic). Additional diagonal regularization is typically only needed near a solution. By contrast, Clarabel uses both static and
dynamic regularization throughout the entire solve to ensure quasi-definiteness and stable factorization. Consequently, Clarabel typically requires at least some iterative refinement to recover the solution of the unregularized matrix, even far away from a solution.

\begin{table}[!htbp]
\centering
\small
\setlength{\tabcolsep}{4pt}
\caption{Median cumulative linear-algebra counts on shared verified feasible
problems. Applications are computed as KKT solves plus refinement steps for
each problem before taking the median; the displayed medians need not add.}
\label{tab:clarabel-linear-work}
\begin{tabular}{@{}llrrrr@{}}
\toprule
Problems & Solver & Factorizations & KKT solves & Refinement steps & Applications \\
\midrule
All ($N=580$) & \cpot{} & 19 & 54 & 1 & 60 \\
 & Clarabel & 15.5 & 45 & 64 & 111.5 \\
\addlinespace
LP ($N=354$) & \cpot{} & 25 & 75 & 2 & 81 \\
 & Clarabel & 20 & 60 & 86 & 147 \\
\addlinespace
QP ($N=226$) & \cpot{} & 11 & 33 & 0 & 35 \\
 & Clarabel & 12 & 35 & 38.5 & 74 \\
\bottomrule
\end{tabular}
\end{table}

\subsubsection{Infeasibility performance}

We compare the solvers on 39 infeasible instances: 29 from NETLIB and ten from
DC-OPF, comprising 37 LPs and two QPs. We vary the absolute and relative
infeasibility tolerances together over $10^{-7},10^{-8},10^{-9},10^{-10}$,
keeping the feasibility and gap tolerances at $10^{-8}$.

We independently verify certificates of primal infeasibility using the original
problem data. To compare certificates independently of their scale, we normalize
each returned certificate with $b^\top y<0$ and require
\[
  \hat y = \frac{y}{-b^\top y},
  \qquad b^\top\hat y=-1,
  \qquad
  \max\!\left\{
    \|A^\top\hat y\|_\infty,
    \|\hat y-\Pi_{\mathcal K^*}(\hat y)\|_\infty
  \right\} \leq 10^{-6}.
\]
Only certificates satisfying this independent criterion count as verified
infeasibility detections.
Table~\ref{tab:clarabel-infeasible} reports the verified certificate counts.
\cpot{} produces independently verified certificates for more instances at
tolerances $10^{-7}$ and $10^{-8}$, and matches Clarabel's counts at
$10^{-9}$ and $10^{-10}$. At tolerance $10^{-8}$, it produces verified
certificates for 36 of the 39 instances, compared with 34 for Clarabel;
32 instances have verified certificates from both solvers.

\begin{table}[!htbp]
\centering
\small
\setlength{\tabcolsep}{4pt}
\caption{Verified primal-infeasibility certificate counts out of 39 instances
at different solver infeasibility tolerances. Normalized certificates are
independently verified at $10^{-6}$.}
\label{tab:clarabel-infeasible}
\begin{tabular}{@{}crr@{}}
\toprule
Infeasibility tolerance & \cpot{} & Clarabel \\
\midrule
$10^{-7}$ & \textbf{36} & 33 \\
$10^{-8}$ & \textbf{36} & 34 \\
$10^{-9}$ & \textbf{35} & \textbf{35} \\
$10^{-10}$ & \textbf{35} & \textbf{35} \\
\bottomrule
\end{tabular}
\end{table}

\begin{figure}[!htbp]
\centering
\includegraphics[width=0.9\linewidth]
{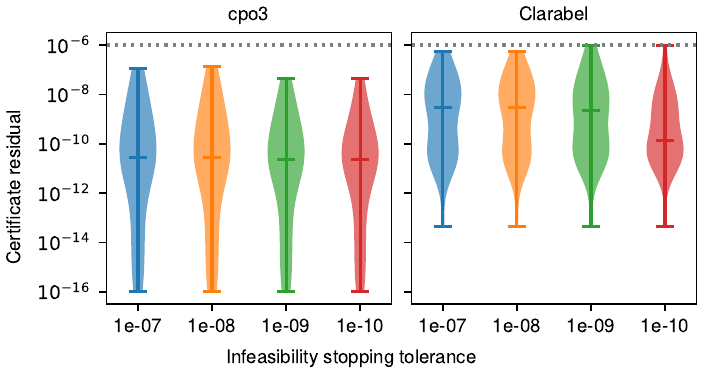}
\caption{Residuals of primal-infeasibility certificates after normalization to
$b^\top\hat y=-1$, at four solver infeasibility tolerances: \cpot{} (left),
Clarabel (right). The dotted line marks the
verification threshold $10^{-6}$, and the interior horizontal marks indicate
medians. Overall \cpot{} achieves certificates with lower residuals than Clarabel.}
\label{fig:clarabel-certificates}
\end{figure}

Figure~\ref{fig:clarabel-certificates} shows the distribution of the normalized
residuals of primal-infeasibility certificates.
Among independently verified certificates, \cpot{} has a lower median normalized residual
at every tested tolerance. At tolerance $10^{-8}$, the median normalized
certificate residual is $2.7\times10^{-11}$ for \cpot{} and
$2.9\times10^{-9}$ for Clarabel, a difference of approximately two orders
of magnitude.

\subsection{Warm-start performance}

The strong monotonicity properties of our central path enabled useful error bounds in Section~\ref{subsec:error_bounds}, which we used to derive a warm start heuristic in Section~\ref{subsec:warm-start}. Here we evaluate
the performance of the warm-start heuristic in two settings.

\subsubsection{Perturbed benchmark problems}
We evaluate this warm-start initialization on perturbed Maros--M\'esz\'aros QPs
and NETLIB LPs. We keep $P$, $A$, and the cones fixed and perturb each entry
of $b$ and $c$ by $\varepsilon_{\mathrm{pert}}$ multiplied by its magnitude, using random signs drawn
once per problem and reused across
$\varepsilon_{\mathrm{pert}}\in\{10^{-4},10^{-3},10^{-2},10^{-1}\}$. After excluding equality-only
problems and unsuccessful baseline solves, the experiment contains 127 QPs
and 91 LPs. Each warm run starts from the original problem's primal--dual
solution, using the construction in Section~\ref{subsec:warm-start}. Table~\ref{tab:warm-start} reports
median iterations and median per-problem speedups
$t_{\mathrm{cold}}/t_{\mathrm{warm}}$ over pairs for which both runs report
\texttt{solved}, with the number of pairs shown explicitly. Some perturbed problems become infeasible or reach the 100-iteration limit, and in 4 of the 872 pairs the cold and warm runs report different statuses. We therefore report iterations and speedups only for pairs where both runs report \texttt{solved}.

\begin{table}[!htbp]
\centering
\small
\caption{Warm-start results for perturbed NETLIB LPs and Maros--M\'esz\'aros
QPs. Pairs is the number of problems for which both runs report \texttt{solved}.
Speedups use solver-reported solve times.}
\label{tab:warm-start}
\begin{tabular}{@{}lcrrrr@{}}
\toprule
Problems & $\varepsilon_{\mathrm{pert}}$ & Pairs & \multicolumn{2}{c}{Median iterations} & Median \\
\cmidrule(lr){4-5}
& & & Cold & Warm & speedup \\
\midrule
LP & $10^{-4}$ & 76 & 25.5 & 8 & $2.29\times$ \\
   & $10^{-3}$ & 75 & 24 & 9 & $1.94\times$ \\
   & $10^{-2}$ & 70 & 24 & 13 & $1.65\times$ \\
   & $10^{-1}$ & 58 & 24 & 20 & $1.21\times$ \\
\addlinespace
QP & $10^{-4}$ & 124 & 18 & 4 & $2.28\times$ \\
   & $10^{-3}$ & 123 & 17 & 6 & $2.03\times$ \\
   & $10^{-2}$ & 117 & 17 & 8 & $1.64\times$ \\
   & $10^{-1}$ & 107 & 17 & 13 & $1.28\times$ \\
\bottomrule
\end{tabular}
\end{table}

Unsurprisingly, warm-starts provide the greatest benefit for small perturbations. At
$\varepsilon_{\mathrm{pert}}=10^{-4}$, median speedups are $2.28\times$ for QPs and $2.29\times$
for LPs, with median iterations falling from 18 to 4 and from 25.5 to 8,
respectively. The benefit decreases as perturbations grow: at
$\varepsilon_{\mathrm{pert}}=10^{-1}$, the median speedup is $1.28\times$ for QPs and $1.21\times$
for LPs.

\subsubsection{Model predictive control (MPC)}
In practice, warm-starts commonly arise when solving a sequence of related problems, such as closed-loop model-predictive control (MPC). Every iteration of MPC solves a QP corresponding to open-loop control of a system, the first action is applied, the resulting next state is observed and the process repeats from the new updated starting point \cite{RawlingsMayneDiehl2017,CamachoBordons2007,Maciejowski2002,BorrelliBemporadMorari2017}. Consequently, each QP is highly related to the previous one. In this example, we control a randomly generated linear system $x_{k+1}=Fx_k+Gu_k$ with $8$ states and $3$ inputs, where $F$ is scaled to spectral radius $0.98$ and the state is driven by Gaussian process noise of standard deviation $0.02$. The MPC problem is
\begin{align}\label{eq:mpc}
    \min_{x,u} \quad & \displaystyle\sum_{k=1}^{N}\|x_k-r\|^2
        +0.1\sum_{k=0}^{N-1}\|u_k\|^2, \nonumber \\
    \text{s.t.} \quad & x_0 = \hat x, \quad  x_{k+1}=Fx_k+Gu_k,\quad k=0,\dots,N-1,\\
    &  \|u_k\|_\infty\leq0.1,\quad
        \|x_{k+1}\|_\infty\leq5,\quad k=0,\dots,N-1, \nonumber\\
    & x_N=x_s, \nonumber
\end{align}
with horizon $N=40$. The current state is denoted by
$\hat x$ and $r$ is a constant setpoint with
$\|r\|_\infty=4$. The setpoint cannot be held within the input bounds, so the terminal state is constrained to the steady state $x_s$ with the smallest stage cost. It minimizes $\|x_s-r\|^2+0.1\|u_s\|^2$ subject to $x_s=Fx_s+Gu_s$ and the same box constraints. Each control step is a QP with $440$ variables, $328$ equality rows and $880$ inequality rows whose data change only through $\hat x$, which enters $b$ in the QP problem data. The setpoint and steady state keep the QPs from becoming trivial as the loop settles, with at least $80$ active inequality constraints for each QP. 

The loop runs for $60$ steps from a state with $\|\hat x\|_\infty = 3.2$, and the QP at every step is solved three ways: cold, warm from the previous step's primal--dual solution, and warm from that solution shifted forward by one stage. The shift drops the first stage and repeats the last for the primal, dual and slack vectors alike. The loop is driven by the cold solve's first input, so all three variants solve the identical sequence of
problems.

\begin{figure}[!htbp]
\centering
\includegraphics[width=\linewidth]{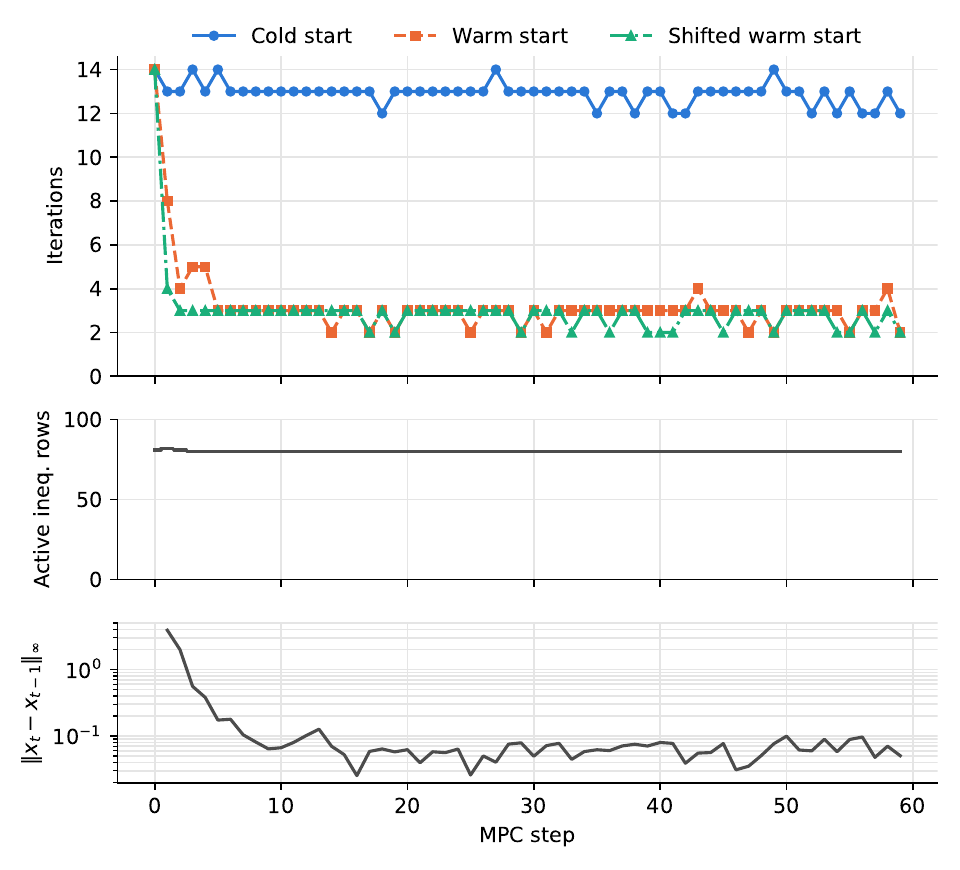}
\caption{Warm starts on a closed-loop run of the MPC
problem~\eqref{eq:mpc} over $60$ steps, each step a QP with $440$
variables, $328$ equality and $880$ inequality rows. Top: iterations per
step for a cold start, and the two warm starts. Middle: the number of inequality rows active
at the cold solution, between $80$ and $82$ input bounds at every step. Bottom: the
change in the state between consecutive steps, $\|x_t - x_{t-1}\|_\infty$,
which is the only change in the problem data from one step to the next.}
\label{fig:warm_mpc}
\end{figure}

Figure~\ref{fig:warm_mpc} shows the iteration count at every step. The cold
start needs a median of $13$ iterations, while the warm and shifted warm
starts take a median of $3$ steps each. Both warm starts are more effective than the cold start, with the shifted warm start being slightly more efficient. Evidently, the warm-start heuristic for our central path is providing a
substantial speedup over cold-starts in this setting, despite the number of active constraints.

\section{Conclusion}

We have introduced a canonical central path for the homogeneous
embedding of quadratic cone programs, defined by strongly monotone
equations with regularization intrinsic to the construction and tied
to the barrier parameter. The path is uniquely defined for every
positive barrier parameter, is smooth, and lies on a sphere of fixed
radius. When a finite primal--dual solution or a strict infeasibility
certificate exists, its limit points as the barrier parameter tends
to zero recover optimal solutions or infeasibility certificates.
Euclidean and local-norm error bounds quantify proximity to the
path and motivate a warm-start heuristic.

The construction permits step computations that retain the nonlinear
embedding equations and linearize only the complementarity conditions
$s\circ y=\mu e$. For linear and quadratic programs, each step reduces
to quasi-definite linear solves and a scalar quadratic equation.
A surrogate formulation allows the predictor and corrector in a
modified Mehrotra scheme to share a single sparse factorization.
The resulting method is implemented in the free and open-source
solver \cpot{}.

Even with presolve and crossover disabled, \cpot{} returns verified
solutions for 631 of 675 instances across eight benchmark collections,
achieving the highest solve count and lowest penalized shifted
geometric mean runtime among the tested solvers. It solves every
Maros--M\'esz\'aros instance and achieves the highest solve counts
and lowest aggregate runtimes on both LPs and QPs.
On problems verified by both \cpot{} and Clarabel, \cpot{}
has lower median counts of refinement corrections and factorization
applications, but a higher median factorization count.
In separate primal-infeasibility experiments, it returns at least as
many verified certificates and lower median certificate residuals
than Clarabel at every tested tolerance. Warm starting yields median
runtime speedups of $2.28\times$ and $2.29\times$ at the smallest
perturbations of Maros--M\'esz\'aros and Netlib problems, respectively,
with gains decreasing as the perturbations grow, and reduces the
iteration count of a closed-loop MPC run with a persistently active set
of input bounds from a median of $13$ to $3$ iterations per step.

These results demonstrate that retaining the nonlinear embedding
structure is compatible with efficient sparse linear algebra and
competitive solver performance. The central-path analysis applies
to general proper cones equipped with logarithmically homogeneous
self-concordant barriers, while the algorithm and implementation
are developed for linear and quadratic programs. Extending the
method to general cones and establishing convergence and complexity
guarantees for the adaptive predictor--corrector scheme remain
important directions for future work.

\bibliographystyle{plain}   
\bibliography{refs}

\clearpage
\appendix
\section{Proofs}\label{app:proofs}

\subsection{Proofs for Section~\ref{sec:central_path_section}}\label{app:path}

\begin{proof}[Proof of Theorem~\ref{thm:central_path_existence}]
The strategy is to realize $T_\mu$ as the restriction of a set-valued operator on $\R^{n+m+1}$ that is maximally
monotone and strongly monotone, whose zero set is then a singleton by standard results.

For this proof, extend $\phi$ to all of $\R^{n+m+1}$ by setting
$\phi(u)=+\infty$ outside $\operatorname{int}(\mathcal C)$.
The extended function is closed, convex and proper, and $\partial\phi$ is
maximally monotone \cite[Thm.~20.25]{Bauschke_Combettes_2017}. Differentiability of $\phi$ on its open domain
gives $\operatorname{dom} \partial \phi = \operatorname{int}( \mathcal{C})$.
Let $\mathcal{Q}$ be the maximal monotone extension of $Q$ from Corollary~\ref{cor:unique_Q}. Recall $\mathcal{Q}(u) = \{ Q(u) \}$ for $u \in \operatorname{int}(\mathcal{C})$ and
$\operatorname{int}(\mathcal{C}) \subset \operatorname{dom}\mathcal{Q}$. Next, define
\begin{equation}\label{eq:max_T_extension}
    \mathcal{T}_\mu \coloneq \mathcal{Q} + \mu \partial \phi + \mu I,
\end{equation}
which is clearly $\mu$-strongly monotone. Noticing that $\operatorname{dom} \mathcal{Q} \cap \operatorname{int}(\operatorname{dom} \partial \phi) = \operatorname{int}(\mathcal{C})$,
which is nonempty, two applications of \cite[Cor.~25.5]{Bauschke_Combettes_2017} deduce that \eqref{eq:max_T_extension} is maximally monotone. We can therefore conclude from
\cite[Cor.~23.37]{Bauschke_Combettes_2017} that $\operatorname{zer}(\mathcal{T}_\mu)$ is a singleton. Furthermore, since
$\operatorname{dom} \partial \phi = \operatorname{int}(\mathcal{C}) \subset
\operatorname{dom} \mathcal{Q}$, we have $\operatorname{dom} \mathcal{T}_\mu = \operatorname{int}(\mathcal{C})$. It therefore must hold that $\mathcal{T}_\mu(u) = \{ Q(u) + \mu \nabla \phi(u) + \mu u \} =
\{ T_\mu(u) \}$. Zeros of $\mathcal{T}_\mu$ are therefore exactly the solutions \eqref{eq:central_path_reduced}, and there is exactly one for any given $\mu > 0$.
\end{proof}

\begin{proof}[Proof of Theorem~\ref{thm:central-path-smoothness}]
The proof relies on the Implicit Function Theorem
\cite[Theorem 1B.1]{dontchev2009implicit} to show local smoothness,
and Theorem~\ref{thm:central_path_existence} to extend the result globally.
Fix $\mu_0>0$. The map $(\mu,u)\mapsto T_\mu(u)$ is of class
$C^d$, and strong monotonicity implies that
$\nabla_u T_{\mu_0}(u_{\mu_0})$ is nonsingular.
The Implicit Function Theorem therefore gives a local solution
map of class $C^d$ near $\mu_0$
\cite[Proposition 1B.5]{dontchev2009implicit}.
By uniqueness, this local map coincides with $\mu\mapsto u_\mu$.
Since $\mu_0$ was arbitrary, $(\mu\mapsto u_\mu)\in C^d(\R_{>0})$.
\end{proof}

\begin{proof}[Proof of Lemma~\ref{lem:limit_points}]
    Fix $\mu_0>0$ and write $R\coloneq\sqrt{\nu+1}$.
    Proposition~\ref{prop:sphere} gives $\|u_\mu\|=R$.
    For $0<\mu\leq\mu_0$, the equation
    \[
        s_\mu=-Ax_\mu+b\tau_\mu+\mu y_\mu
    \]
    therefore bounds $s_\mu$ uniformly. The final block
    equation gives
    \[
        \kappa_\mu+
        \frac{x_\mu^\top Px_\mu}{\tau_\mu}
        =
        -c^\top x_\mu-b^\top y_\mu+\mu\tau_\mu.
    \]
    Since both terms on the left are nonnegative,
    \[
        0\leq\kappa_\mu
        \leq R\bigl(\|(c,b)\|+\mu_0\bigr).
    \]
    Thus $\{(u_\mu,v_\mu):0<\mu\leq\mu_0\}$ is bounded,
    so it has a limit point as $\mu\downarrow0$.

    Let $(u_{\mu_j},v_{\mu_j})\to(u,v)$ along
    a sequence $\mu_j\downarrow0$. Since
    \[
        v_{\mu_j}-\mu_j u_{\mu_j}
        \in\mathcal Q(u_{\mu_j}),
    \]
    the closed graph of $\mathcal Q$
    \cite[Prop.~20.38]{Bauschke_Combettes_2017}
    gives $v\in\mathcal Q(u)$.
    Closedness of the cones gives
    $u\in\mathcal C$ and $v\in\mathcal C^*$.
    Finally, logarithmic homogeneity yields
    \[
        u_{\mu_j}^\top v_{\mu_j}
        =
        -\mu_j u_{\mu_j}^\top\nabla\phi(u_{\mu_j})
        =
        \mu_j(\nu+1)\to0.
    \]
    Hence $u^\top v=0$, and
    $(u,v)$ solves
    \eqref{eq:embedding_compact}.
\end{proof}

\begin{proof}[Proof of Theorem~\ref{thm:barrier_error_bound}]
    Since $\nabla^2\Psi(u)=I+\nabla^2\phi(u)\succeq I$, the Hessian is
    nonsingular on $\operatorname{int}(\mathcal{C})$, so $\eta_\mu(u)$ is well
    defined and $\|\cdot\|_{\nabla^2\Psi(u)}$ is a norm dominating $\|\cdot\|$,
    which proves the first inequality. Let $r\coloneq u-u_\mu$ and $t\coloneq\|r\|_{\nabla^2\Psi(u)}$. The case $t = 0$ is immediate. Since both
    $u$ and $u_\mu$ lie in $\operatorname{int}(\mathcal{C})$,
    \cite[Thm.~4.1.7]{nesterov2004lectures} gives
    \begin{equation}\label{eq:nesterov_417}
        r^\top\bigl(\nabla\Psi(u)-\nabla\Psi(u_\mu)\bigr)\ge\frac{t^2}{1+t}.
    \end{equation}
    Because $T_\mu=Q+\mu\nabla\Psi$ and $T_\mu(u_\mu)=0$,
    \[
        r^\top T_\mu(u)
        =r^\top\bigl(Q(u)-Q(u_\mu)\bigr)
        +\mu\,r^\top\bigl(\nabla\Psi(u)-\nabla\Psi(u_\mu)\bigr)
        \ge\mu\,\frac{t^2}{1+t},
    \]
    where the first term is nonnegative by monotonicity of $Q$ and
    the second is bounded by \eqref{eq:nesterov_417}. On the other hand,
    Cauchy--Schwarz in the metric $\nabla^2\Psi(u)$ gives
    \[
        r^\top T_\mu(u)
        \le\|r\|_{\nabla^2\Psi(u)}\,\|T_\mu(u)\|_{(\nabla^2\Psi(u))^{-1}}
        =\mu\,t\,\eta_\mu(u).
    \]
    Combining the two bounds and dividing by $\mu t>0$ yields
    $t/(1+t)\le\eta_\mu(u)$, and since $\eta_\mu(u)<1$ this rearranges to
    $t\le\eta_\mu(u)/(1-\eta_\mu(u))$.
\end{proof}

\begin{proof}[Proof of Theorem~\ref{thm:nondegenerate_limits}]
    Recall that
    \[
        Q(u_\mu)=v_\mu-\mu u_\mu,\qquad
        u_\mu^\top v_\mu=\mu(\nu+1),\qquad
        \|u_\mu\|^2=\nu+1,
    \]
    and $\tau_\mu\kappa_\mu=\mu$.
    For any complementary embedding pair $(u,v)$ with
    $u\in\mathcal C$, $v\in\mathcal C^*$,
    $v\in\mathcal Q(u)$ and $u^\top v=0$,
    monotonicity of $\mathcal Q$ gives
    \[
        0
        \leq (u_\mu-u)^\top(v_\mu-\mu u_\mu-v)
        =\mu u^\top u_\mu-u_\mu^\top v-u^\top v_\mu.
    \]
    Both $u_\mu^\top v$ and $u^\top v_\mu$ are
    nonnegative by cone duality. Rearranging and applying
    Cauchy--Schwarz therefore yields
    \begin{equation}\label{eq:embedding_comparison}
        0
        \leq u_\mu^\top v+u^\top v_\mu
        \leq \mu u^\top u_\mu
        \leq \mu\sqrt{\nu+1}\|u\|.
    \end{equation}

    For \textit{(i)}, take $u=(x,y,1)$ and
    $v=(0,s,0)$.
    The KKT conditions imply that $(u,v)$ is a
    complementary embedding pair, with $Q(u)=v$.
    Since
    \[
        u^\top v_\mu
        =y^\top s_\mu+\kappa_\mu
        \geq\kappa_\mu
        =\frac{\mu}{\tau_\mu},
    \]
    \eqref{eq:embedding_comparison} gives
    \[
        \frac{\mu}{\tau_\mu}
        \leq
        \mu\sqrt{\nu+1}
        \sqrt{1+\|x\|^2+\|y\|^2}.
    \]
    Dividing by $\mu>0$ and rearranging proves the
    lower bound on $\tau_\mu$.

    For \textit{(ii)}, take $u=(0,y,0)$ and
    $v=(0,0,-b^\top y)$.
    Since $y\in\mathcal P$, we have
    $u\in\mathcal C$, $v\in\mathcal C^*$ and
    $u^\top v=0$. Moreover, $A^\top y=0$ gives
    \[
        Q(0,y,t)\to v
        \qquad\text{as }t\downarrow0.
    \]
    The graph of $\mathcal Q$ is closed
    \cite[Prop.~20.38]{Bauschke_Combettes_2017}, so
    $v\in\mathcal Q(u)$.
    Applying \eqref{eq:embedding_comparison} yields
    \[
        (-b^\top y)\tau_\mu
        =u_\mu^\top v
        \leq\mu\sqrt{\nu+1}\|y\|.
    \]
    Dividing by $-b^\top y>0$ and using
    $\kappa_\mu=\mu/\tau_\mu$ gives
    \[
        \tau_\mu
        \leq\frac{\mu\sqrt{\nu+1}\|y\|}{-b^\top y},
        \qquad
        \kappa_\mu
        \geq\frac{-b^\top y}{\sqrt{\nu+1}\|y\|}.
    \]

    For \textit{(iii)}, take $u=(x,0,0)$ and
    $v=(0,-Ax,-c^\top x)$.
    Since $x\in\mathcal D$, we have
    $u\in\mathcal C$, $v\in\mathcal C^*$ and
    $u^\top v=0$. Moreover, $Px=0$ gives
    \[
        Q(x,0,t)\to v
        \qquad\text{as }t\downarrow0,
    \]
    so closedness of the graph again implies
    $v\in\mathcal Q(u)$.
    Cone duality gives
    \[
        u_\mu^\top v
        =y_\mu^\top(-Ax)+(-c^\top x)\tau_\mu
        \geq(-c^\top x)\tau_\mu.
    \]
    Hence \eqref{eq:embedding_comparison} yields
    \[
        (-c^\top x)\tau_\mu
        \leq\mu\sqrt{\nu+1}\|x\|,
    \]
    and therefore
    \[
        \tau_\mu
        \leq\frac{\mu\sqrt{\nu+1}\|x\|}{-c^\top x},
        \qquad
        \kappa_\mu
        \geq\frac{-c^\top x}{\sqrt{\nu+1}\|x\|}.
    \]

    Finally, Lemma~\ref{lem:limit_points} guarantees the
    existence of limit points as $\mu\downarrow0$ and
    shows that every such point solves
    the homogeneous embedding. In case \textit{(i)},
    the uniform lower bound on $\tau_\mu$ gives a positive
    $\tau$ component at every limit point. In cases \textit{(ii)}
    and \textit{(iii)}, the uniform lower bound on $\kappa_\mu$
    gives a positive $\kappa$ component.
    Lemma~\ref{lem:asymp_feasible} therefore recovers a
    primal--dual optimal solution in case \textit{(i)}
    and an infeasibility certificate in cases
    \textit{(ii)} and \textit{(iii)}.
\end{proof}

\subsection{Proofs for Section~\ref{sec:ipm_section}}\label{app:ipm}

\begin{proof}[Proof of Proposition~\ref{thm:resolvent}]
Write $K$ for the coefficient matrix of \eqref{eq:qd_system}.
Since $P+\rho I\succ0$ and $D_k+\rho I\succ0$, $K$ is
symmetric quasi-definite and hence nonsingular
\cite{vanderbei1995}, so the four vectors in
\eqref{eq:qd_system} are well defined.

First suppose $\tau_{k+1}>0$. Negating the second block row,
the first two rows of \eqref{eq:generic_system} read
\[
    K
    \begin{bmatrix} x_{k+1}\\ y_{k+1} \end{bmatrix}
    =
    \begin{bmatrix} r_x\\ -r_y \end{bmatrix}
    +\tau_{k+1}
    \begin{bmatrix} -c\\ b \end{bmatrix}.
\]
By invertibility of $K$, this holds if and only if
\[
    x_{k+1}=x'+\tau_{k+1}x'',
    \qquad
    y_{k+1}=y'+\tau_{k+1}y''.
\]
Substituting into the third row of \eqref{eq:generic_system}
and multiplying by $\tau_{k+1}>0$ gives
\begin{equation}\label{eq:tau_quadratic}
    \theta_2\tau_{k+1}^2+\theta_1\tau_{k+1}-\theta_0=0.
\end{equation}
Thus solutions of \eqref{eq:generic_system} with
$\tau_{k+1}>0$ correspond bijectively to positive roots of
\eqref{eq:tau_quadratic}.

Multiplying the first block row of the $(x'',y'')$ column of
\eqref{eq:qd_system} by $(x'')^\top$, the second by $(y'')^\top$,
and subtracting gives
\[
    -c^\top x''-b^\top y''
    =
    (x'')^\top(P+\rho I)x''
    +(y'')^\top(D_k+\rho I)y''.
\]
Consequently,
\[
    \theta_2
    =
    \rho\bigl(1+\|x''\|^2+\|y''\|^2\bigr)
    +(y'')^\top D_k y''
    \geq\rho>0.
\]
Also, $\theta_0=(x')^\top Px'+\mu\geq\mu\geq0$.

If $\theta_0>0$, the roots of \eqref{eq:tau_quadratic}
have product $-\theta_0/\theta_2<0$, so exactly one is
positive, and it is selected by the stated quadratic formula.
If $\theta_0=0$ and $\theta_1<0$, the roots are $0$ and
$-\theta_1/\theta_2>0$, and the formula returns the latter.
In either case, the recovered triple is the unique solution
of \eqref{eq:generic_system} with $\tau_{k+1}>0$.

It remains to consider $\theta_0=0$ and $\theta_1\geq0$.
The roots are then $0$ and $-\theta_1/\theta_2\leq0$, so
there is no positive solution and the formula returns
$\tau_{k+1}=0$. Since
\[
    (x')^\top Px'+\mu=0,
    \qquad P\succeq0,\quad \mu\geq0,
\]
we have $\mu=0$ and $Px'=0$. Hence
$u_{k+1}=(x',y',0)$ lies in the domain of the boundary
operator $\mathcal I$, and
\[
    \theta_1=-c^\top x'-b^\top y'-r_\tau\geq0.
\]
Together with the first two block equations in
\eqref{eq:qd_system}, this gives
\[
    r-\rho u_{k+1}-(0,D_k y',0)
    =
    (A^\top y',-Ax',r_\tau)
    \in\mathcal I(x',y',0)
    =\mathcal Q(u_{k+1}),
\]
where the inclusion follows from \eqref{eq:I_operator}.
Thus $u_{k+1}$ solves the extended system
\[
    r\in\mathcal Q(u_{k+1})
        +\rho u_{k+1}+(0,D_k y_{k+1},0).
\]
Since $\mathcal Q$ is maximal monotone and $D_k\succeq0$, the
operator on the right is $\rho$-strongly monotone, so this
solution is unique.
\end{proof}

\begin{proposition}[Contraction of the surrogate system]
    \label{prop:contraction}
    Let $\mu_k>0$, $\mu\in[0,\mu_k]$, and fix the same
    offset $\delta$ in \eqref{eq:split_system} and
    \eqref{eq:surrogate_system}.
    Suppose $u^\star$ solves \eqref{eq:split_system}
    and $u_{k+1}$ solves \eqref{eq:surrogate_system}
    for a surrogate $\tilde u$, with both points in
    $\operatorname{dom}Q$. Then
    \[
        \|u_{k+1}-u^\star\|
        \leq
        \left(1-\frac{\mu}{\mu_k}\right)
        \|\tilde u-u^\star\|.
    \]
\end{proposition}

\begin{proof}[Proof of Proposition~\ref{prop:contraction}]
The subproblem residual $S_k$ is $\mu$-strongly monotone.
Subtracting $S_k(u^\star)=0$ from the surrogate equation gives
\[
    S_k(u_{k+1})-S_k(u^\star)
    +(\mu_k-\mu)(u_{k+1}-u^\star)
    =(\mu_k-\mu)(\tilde u-u^\star).
\]
Taking the inner product with $u_{k+1}-u^\star$, applying
strong monotonicity on the left and Cauchy--Schwarz on the
right yields
\[
    \mu_k\|u_{k+1}-u^\star\|^2
    \leq(\mu_k-\mu)\|\tilde u-u^\star\|\,\|u_{k+1}-u^\star\|.
\]
If $u_{k+1}=u^\star$, there is nothing to prove. Otherwise,
dividing by $\mu_k\|u_{k+1}-u^\star\|$ gives
\[
    \|u_{k+1}-u^\star\|
    \leq\left(1-\frac{\mu}{\mu_k}\right)\|\tilde u-u^\star\|.
    \qedhere
\]
\end{proof}

\clearpage
\section{Presolve}
\label{app:presolve}

The following LP reductions from PSLP~\cite{cederberg2026presolving}
apply unchanged to quadratic programs because they depend only
on the constraints and variable bounds, not on the objective:
singleton rows, redundant constraints, parallel rows,
and primal propagation.

For variables absent from the quadratic part of the objective
($P_{\cdot j}=0$),
PSLP's rules for column singletons in inequalities and variable
locks also apply unchanged. Below, we describe how
other reductions account for quadratic terms.

\paragraph{Variable substitutions.}
We eliminate a variable from an equality constraint if it is one of
only two variables in that equality, or if it appears in no other constraint.
We solve the equality for the selected variable and substitute
the resulting expression into the remaining constraints and objective.

For a QP, this substitution also updates the quadratic terms,
including couplings between variables. Elimination may introduce
new nonzero coefficients in the constraints and objective, but
reduces the overall number of decision variables.

\paragraph{Variable-lock fixing.}
Let $\ell \leq x \leq u$ denote the current variable bounds,
whose lower and upper endpoints may be $-\infty$ and
$+\infty$, respectively. A variable $x_j$ may be fixed at a finite
lower bound $\ell_j$ if decreasing
$x_j$ violates no other constraint and cannot increase the objective.
For an LP, the objective condition is $c_j \geq 0$.
For a QP, we compute a lower bound on the derivative with respect to $x_j$
at $x_j=\ell_j$ by minimizing each contribution from the other variables:
\[
    P_{jj}\ell_j+c_j+
    \sum_{k \ne j}
    \min\{P_{jk}\ell_k, P_{jk}u_k\}.
\]
We use the convention that $0 \cdot \pm \infty = 0$.
If this bound is nonnegative, convexity implies that the
objective is nondecreasing in $x_j$ throughout its allowed
interval, for every admissible choice of the other variables.
Thus, decreasing $x_j$ to $\ell_j$ cannot increase the objective.
The upper-bound rule is analogous: increasing $x_j$ must
violate no other constraint, and the derivative's upper bound
at $x_j=u_j$, obtained by replacing $\min$ with $\max$,
must be nonpositive.

\paragraph{Parallel columns.}
For proportional constraint columns $A_{\cdot k}=\omega A_{\cdot j}$,
we use parallel-column aggregation~\cite[Section~6.3]{achterberg2020presolve}
and the proportional-column special case of the bound-fixing
rule in~\cite[Section~6.4]{achterberg2020presolve}.
For a QP, we additionally require $P_{\cdot k}=\omega P_{\cdot j}$.
This ensures that increasing $x_k$ by $t$ and decreasing $x_j$
by $\omega t$ leaves the quadratic term unchanged, so the
LP exchange argument remains valid.

\paragraph{Variables subject only to bounds.}
A variable that appears in no constraint except simple bounds
may still be coupled to other variables through the quadratic objective.
If these couplings are absent and $P_{jj}>0$, we minimize
its scalar quadratic objective by fixing $x_j$ at
$-c_j/P_{jj}$, clipped to its bounds.

If quadratic couplings are present, we can eliminate $x_j$
when it has no bounds and $P_{jj}>0$. Setting its objective
derivative to zero gives
\[
    x_j=-\frac{c_j+\sum_{k\ne j}P_{jk}x_k}{P_{jj}}.
\]
Substituting this expression into the objective eliminates
$x_j$ and replaces the remaining Hessian by its Schur complement.
This elimination may introduce fill-in in the remaining Hessian,
so we apply it only when the total number of Hessian nonzeros
does not increase.

\clearpage
\section{Performance with default presolve and crossover settings}
\label{app:default-settings}
\setlength{\intextsep}{8pt}

We repeat the comparison of Section~\ref{sec:central-path}
using the solvers' default presolve and crossover settings
where available. All other experimental settings and
performance measures are unchanged.

\begin{figure}[!htbp]
\centering
\includegraphics[width=0.9\linewidth]
{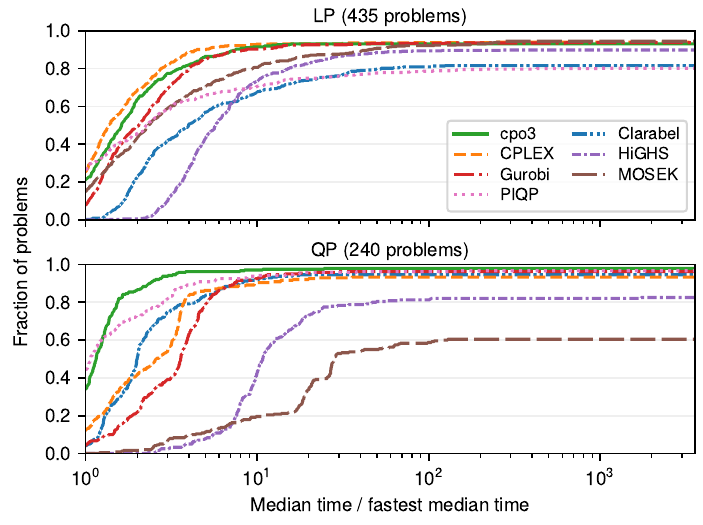}
\caption{Dolan--Mor\'e performance profiles with default presolve and crossover settings:
435 LPs (top); 240 QPs (bottom). The solver legend applies to both panels.}
\label{fig:default-settings}
\end{figure}

\cpot{} has the most verified solves (640/675) among all tested solvers and the
lowest penalized shifted geometric mean runtime among open-source alternatives.
On QPs, it achieves both the highest verified solve count (235/240) and the
lowest mean runtime (1.28 seconds). On LPs, it verifies 405 of 435 problems,
six fewer than MOSEK's leading total. Table~\ref{tab:default-settings}
summarizes the results, and Figure~\ref{fig:default-settings} shows the
performance profiles.

\begin{table}[!htbp]
\centering
\small
\caption{Verified solves and shifted geometric mean runtimes (10-second shift,
300-second failure penalty). Best results in each column are shown in bold.}
\label{tab:default-settings}
\begin{tabular}{@{}lrrrrrr@{}}
\toprule
& \multicolumn{2}{c}{All ($N=675$)}
& \multicolumn{2}{c}{LP ($N=435$)}
& \multicolumn{2}{c}{QP ($N=240$)} \\
\cmidrule(lr){2-3}\cmidrule(lr){4-5}\cmidrule(lr){6-7}
Solver & Solved & Time (s) & Solved & Time (s) & Solved & Time (s) \\
\midrule
\cpot{} & \textbf{640} & 4.40 & 405 & 6.48 & \textbf{235} & \textbf{1.28} \\
CPLEX & 633 & 4.52 & 409 & 5.18 & 224 & 3.41 \\
Gurobi & 639 & \textbf{3.88} & 408 & \textbf{4.97} & 231 & 2.10 \\
PIQP & 582 & 8.88 & 350 & 13.96 & 232 & 2.25 \\
Clarabel & 582 & 10.32 & 355 & 16.08 & 227 & 2.92 \\
HiGHS & 589 & 9.74 & 391 & 9.88 & 198 & 9.47 \\
MOSEK & 556 & 11.38 & \textbf{411} & 5.30 & 145 & 29.20 \\
\bottomrule
\end{tabular}
\end{table}

\setlength{\intextsep}{4pt}
\setlength{\floatsep}{6pt}

\end{document}